\documentclass[11pt]{article}
\usepackage{graphics}
\usepackage{color}
\usepackage[all]{xy}

\newfont{\frak}{eufm10 scaled\magstep1}

\newfont{\sfrak}{eufm8 scaled\magstephalf}
\newfont{\bbb}{msbm10 scaled\magstephalf}
\newfont{\sbbb}{msbm7 scaled\magstephalf}

\def\cdd{\C^d_{\Delta}}
\def\pic{\pi_{\SC}}
\def\Pic{\Pi_{\SC}}

\def\squo{\Psi^{-1}(0)/N}
\def\cquo{\cdd//\Nc}
\def\zset{\Psi^{-1}(0)}

\def\Dc{D_{\SC}}
\def\dc{\d_{\SC}}

\def\Tdc{T^d_{\SC}}

\def\Nc{N_{\SC}}

\def\D{\Delta}
\def\C{\mbox{\bbb{C}}}
\def\R{\mbox{\bbb{R}}}
\def\Z{\mbox{\bbb{Z}}}
\def\K{\mbox{\bbb{K}}}
\def\Q{\mbox{\bbb{Q}}}

\def\rd{\R^d}

\def\vz{\underline{z}}
\def\vx{\underline{x}}

\def\vu{\underline{u}}
\def\vv{\underline{v}}
\def\vw{\underline{w}}

\def\ed{e_1,\ldots,e_d}

\def\xd{X_1,\ldots,X_d}
\def\d{\mbox{\frak d}}
\def\n{\mbox{\frak n}}
\def\ddu{\d^*}

\def\zd{(z_1,\cdots,z_d)}
\def\lorw{\longrightarrow}
\def\SC{\mbox{\sbbb{C}}}
\def\Dc{D_{\SC}}

\def\zset{\Psi^{-1}(0)}

\def\s{\mbox{\frak{s}}}
\def\r{\mbox{\frak{r}}}

\def\b{\mbox{\frak{b}}}

\def\vt{\tilde{V}}

\def\va{V_{\alpha}}
\def\vta{\tilde{V}_{\alpha}}

\def\ga{\Gamma_{\alpha}}

\def\ga{\Gamma_{\alpha}}

\def\G{\Gamma}

\def\ft{\tilde{f}}

\def\s{{\mbox{\frak{s}}}}

\def\Ker{\hbox{Ker}}
\def\Im{\hbox{Im}}

\def\squareforqed{\hbox{\rlap{$\sqcap$}$\sqcup$}}
\def\qed{\ifmmode\else\unskip\quad\fi\squareforqed}
\def\smartqed{\def\qed{\ifmmode\squareforqed\else{\unskip\nobreak\hfil
\penalty50\hskip1em\null\nobreak\hfil\squareforqed
\parfillskip=0pt\finalhyphendemerits=0\endgraf}\fi}}

\newtheorem{thm}{Theorem}[section]
\newtheorem{prop}[thm]{Proposition}
\newtheorem{lemma}[thm]{Lemma}
\newtheorem{defn}[thm]{Definition}
\newtheorem{remark}[thm]{Remark}
\newtheorem{example}[thm]{Example}
\newtheorem{cor}[thm]{Corollary}
\newcommand{\proof}{\mbox{\textbf{ Proof.\ \ }}}

\title{\sc Geometric spaces from arbitrary convex polytopes}
\author{\sc Fiammetta Battaglia}
\date{23/04/10}
\date{}
\begin{document}
\maketitle

\begin{abstract}
\let\thefootnote\relax\footnotetext{Research partially supported by MIUR (``Geometria Differenziale e Analisi
Globale'' PRIN 2007) and INdAM (GNSAGA)}
We associate a geometric space to an arbitrary convex polytope. Our 
construction parallels the construction by D. Cox of a toric variety as a GIT
quotient \cite{cox}. The spaces that we obtain are endowed with a natural 
stratification and perfectly mimic the features of toric varieties associated 
to rational convex polytopes. 
\end{abstract}

{\parindent 0pt

{\small 2000 \textit{Mathematics Subject Classification.} Primary: 14M25.
Secondary: 32S99.}


 {\small \textit{Affiliation}: Dipartimento di Matematica Applicata, 
Via S. Marta 3, 50139 Firenze, Italy.

\textit{e-mail}: fiammetta.battaglia@unifi.it}
}

\section*{Introduction}
The goal of this paper is to construct geometric spaces from nonrational 
nonsimple polytopes that generalize toric varieties.
To each $n$-dimensional convex polytope $\D$ there corresponds a unique fan, generated by the $1$-dimensional
cones dual to the facets of the polytope. The polytope $\D$ is 
{\em rational} if these $1$-dimensional cones are generated by vectors
lying in a lattice $L$.
The polytope $\D$, in a neighboorhood of a $p$-dimensional face $F$, 
is the product of $F$ by a cone over an $(n-p-1)$-dimensional polytope $\D_F$. 
The face $F$ is regular if $\D_F$ is a simplex, singular otherwise.
The polytope $\D$ is {\em simple} if all of its vertices are regular, which
implies that all of its faces are also regular. 
A simple nonrational convex polytope can always be perturbed into a rational 
one, combinatorially equivalent. This is not true in the nonsimple case:
there are nonsimple convex polytopes that are not even combinatorially equivalent 
to rational ones (the first example, due to M. Perles, was published in the book by B. Gr\"umbaum 
\cite{grun}, see also \cite{ziegler} and \cite{gebert}).

A rational convex polytope in a lattice $L$ gives rise to  
a toric variety $X$, acted on by the torus
$\d_{\SC}/L\simeq(\C^*)^n$, where $\d=L\otimes_{\Z}\R$. There is a one-to-one
correspondence between $p$-dimensional orbits of the torus and $p$-dimensional faces 
of the polytope. The orbit corresponding to the interior of the polytope, $\d_{\SC}/L\simeq(\C^*)^n$, 
is open and dense; its compactification $X$ is obtained by gluing, to the $n-$dimensional  
orbit, the smaller orbits corresponding to the other faces of the polytope. 
The gluing pattern reflects the combinatorics of the polytope.
This subdivision in orbits yields a natural stratification of $X$:
the union of the orbits corresponding to the regular faces is the set
of rationally smooth points of $X$, whilst the orbits corresponding to the
singular faces are the singular strata.

We are interested in the following problem: is there an analogue of this 
construction for nonrational polytopes and what are the spaces we end up with?

When dealing with nonrational polytopes, the first step is to replace lattices
with quasilattices, this was first done by E. Prato in \cite{p}, 
in order to generalize to nonrational simple convex polytopes the Delzant procedure \cite{delzant}. 
A quasilattice $Q$ in a vector space $\d$ is a $\Z$-submodule of 
$\d$ generated by a finite set of spanning vectors of $\d$, it is an important notion
in the physics of quasicrystals (\cite{S},\cite{steurer}). We will also need the notions
of {\em quasitorus} and {\em quasifold}, introduced in \cite{p}. 
Quasifold is a generalization of manifold and orbifold:
the local models are quotients of manifolds modulo the smooth action of discrete groups, not necessarily finite. 
Therefore a quasifold might be non Hausdorff; by compact space we shall always mean a space such that each open covering 
admits a finite  subcovering. Quasitorus is a natural generalization of torus, an example, using the above notations, 
is $\d/Q$.

Let us now go back to the problem posed above: let $\d$ be an $n$-dimensional real vector space, 
we consider a dimension $n$, nonrational, nonsimple convex polytope $\D\subset\d^*$,
together with the following data: a quasilattice $Q$ in $\d$ and a set of generators of the 
$1$-dimensional dual cones, contained in $Q$ -- in the rational case, 
given $\D$ in $L$, the generators are chosen to be primitive in $L$. 
We then consider the complex quasitorus $\d_{\SC}/Q$.
We prove that, to each choice of generators and quasilattice $Q$, there corresponds
a space $X_{\D}$, acted on continuosly by the quasitorus  $\d_{\SC}/Q$.
As in the rational case there is a one-to-one correspondence between
$p$-dimensional orbits of the quasitorus and $p$-dimensional faces of the polytope.
The space $X_{\D}$ is the compactification of the $n$-dimensional
complex quasitorus $\d_{\SC}/Q$, obtained by gluing, to
the open and dense orbit corresponding to the interior of the polytope,
the smaller orbits corresponding to the other faces of the polytope.
This defines a stratification of $X_{\D}$:
the orbits corresponding to the singular faces are the singular strata, 
they are complex quasifolds, isomorphic to $(\C^*)^p$ modulo the action of a discrete 
group; the maximal stratum, given by the union of the orbits 
corresponding to the regular faces, has the structure of an $n$-dimensional 
complex quasifold.
We obtain in fact a {\em complex stratification by quasifolds}, 
that is, in a neighborhood of the stratum corresponding to a singular face $F$, 
the space $X_{\D}$ is biholomorphic to the twisted product,
under the action of a discrete group, of the stratum itself
by a {\em complex cone} over the space $X_{\D_F}$;
when the polytope is rational we recover the stratified structure
of the toric variety, since the twisting group is in this case finite.
We therefore observe on the space $X_{\D}$ two different kinds of 
singularities: the stratified structure of $X_{\D}$ and 
the quasifold structure of the strata, due to nonsimpleness and nonrationality
respectively. Let us see how these two features of $\D$ intervene in
the quotient construction that produces our space $X_{\D}$.

In \cite{cox} D. Cox constructs the toric variety corresponding to a rational
convex polytope $\D$ in a lattice $L$ as the categorical quotient of a suitable 
open subset $C^d_{\D}$ of $\C^d$, modulo the action of a subtorus of $(\C^*)^d$, 
where $d$ is the number of facets of the polytope, or, in other terms, the number 
of $1$-dimensional cones in the dual fan. 
The open subset $\C^d_{\D}$ can be defined, in purely combinatorial terms, 
for any convex polytope. In order to construct, in our setting, a 
suitable subgroup $N_{\SC}$ of $\C^d$, we adopt the generalization of the Delzant procedure 
given in \cite{p}.
The group $N_{\SC}$ thus obtained is the complexification of a nonclosed 
subgroup $N$ of $(S^1)^d$, namely $N_{\SC}=\exp(i\n)N$, where 
$\n=\hbox{Lie}(N)$. 
In the simple case the orbits of $\exp(i\n)$ are closed; 
the geometric quotient, $\C^d_{\D}/N_{\SC}$, gives a rationally 
smooth toric variety 
in the rational case and an $n$--dimensional compact complex quasifold 
in the nonrational case -- this was proved, jointly with Elisa Prato, in \cite{cx}.
But in the nonsimple case there are nonclosed $\exp(i\n)$-orbits. To controle 
their behavior we make use of 
suitable functions on $\C^d$ of the kind 
$$\sum_{k=1}^{d} |z_k|^{c_k},$$ $c_k\in\R$.
These functions turn out to be an essential tool and play  
the role of the rational functions used by Cox in his paper.
We are then able to generalize
the notion of categorical quotient: the space $X_{\D}$ is finally defined to be
the quotient $\C^d_{\D}//N_{\SC}$. 
Remark that there are two distinct ways in which 
$N_{\SC}$-orbits are nonclosed, on one hand $N$ itself is a nonclosed subgroup 
of $(S^1)^d$; on the other hand, since the  polytope is nonsimple, there are nonclosed 
$\exp(i\n)$-orbits: the two kind of singularities of  $X_{\D}$ described above -- quasifold structure
of strata and stratification --
are a direct consequence of this fact.

We continue by considering the symplectic quotient corresponding to $\D$,
with the same choice of generators and quasilattice. This is  a space 
stratified by symplectic quasifolds \cite{ns}. We prove that $X_{\D}$ is 
homeomorphic to its symplectic counterpart, $M_{\D}$, that such homeomorphism 
is compatible with the stratifications of $X_{\D}$ and $M_{\D}$ and its 
restriction to each stratum is a diffeomorphism, with respect to which the 
symplectic and complex structures 
of strata are compatible. In particular the space $X_{\D}$ is compact  
and its strata are in fact K\"ahler quasifolds. 

In conclusion these results, which were announced in \cite{cr}, 
enable us to associate
a geometric space to an arbitrary convex polytope, generalizing the 
construction of a toric variety from a rational polytope;
indeed, in the case of a rational polytope $\D$ in a lattice $L$, 
with the choice of primitive generators, our space $X_{\D}$ 
{\em is} the toric variety associated to the pair $(\D,L)$.  
The geometry and topology of our spaces and the connection with the 
combinatorics of the associated polytopes are natural questions related 
to our work. An initial step towards a better understanding of these different 
aspects, that we are pursuing, is the study of cohomological invariants of 
our spaces. A first result in this direction can be found in \cite{betti}.

\vspace{.3cm}

{\bf Acknowledgments}. We would like to thank the referee of the research 
announcement relative to the present paper, \cite{cr}, for her/his helpful 
remarks.  

\vspace{.3cm}

\section{Singularity types}
Before going into the core of the paper, we recall in this section the 
necessary notions, they can all be found in detail in the references. 
\subsection{Quasifolds}
For the detailed definitions of real quasifolds and related geometrical objects 
we refer the reader to the original article by E. Prato \cite{p} and to the Appendix of the 
joint paper with E. Prato \cite{kite}, where some of the definitions were reformulated. 
The definition of complex quasifold can be derived naturally from that of real quasifold, 
a version based on \cite{p} was given jointly with E. Prato in \cite{cx}.
The definitions of quasitori and their Hamiltonian actions were introduced in \cite{p}, 
an extension of these notions to the complex set up was then given in \cite{cx}. 

Let us briefly recall the definition of complex quasifold. 
Let $\vt$ be a complex quasifold and let $\G$ be a discrete group acting
on $\vt$ by biholomorphims in such a way that
the set of points where the action is not free, is closed
and has minimal real codimension $\geq2$. This implies that the action is free on an open, dense and
connected set. The quotient space $\vt/\G$ is a {\em quasifold model}.
Two models $\vt/\G$ and $\tilde{W}/\D$ are {\em biholomorphic} if
there exists a homeomorphism $f\,\colon\,\vt/\G\rightarrow\tilde{W}/\D$ that lifts
to a biholomorphism $\ft\,\colon\,\vt\rightarrow\tilde{W}$. 
Let $X$ be a topological space, a complex {\em quasifold chart} on $X$ is an open subset
$V\subset X$, homeomorphic to a model $\vt/\G$. Given a pair of intersecting charts a suitable notion
of holomorphic {\em change of charts} is defined, if satisfied the two charts are said to be {\em compatible}.
A topological space is endowed with the structure of a {\em complex quasifold} if it is covered
by a collection ${\mathcal A}=\{\va\simeq\vta/\ga\;|\;\alpha\in A\}$ of compatible complex charts.

Geometric objects on quasifolds, like differential forms, are defined 
on each $\vta$ with the additional conditions that they 
descend to the quotient $\vta/\ga$ and that, when charts are overlapping,
they glue suitably by means of the changes of charts.

As mentioned in the introduction quasilattices are very important: 
\begin{defn}[Quasilattice] Let $\d$ be a real vector space of dimension $n$.
A {\em quasilattice} $Q$ in $\d$ is a $\Z$-submodule of $\d$, 
generated by a set of generators of $\d$.
\end{defn}
Notice that if $\hbox{rank}(Q)=n$ then $Q$ is a lattice.
\begin{defn}[Quasitorus, quasi-Lie algebra, exponential\cite{p,cx}]\label{quasitorus}
{\rm Let $\d$ be a vector space of dimension $n$ and let $\dc=\d+i\d$ its 
complexification. Let $Q$ be a quasilattice in $\d$. The quotient
$D=\d/Q$ (respect. $\Dc=\dc/Q$) is an $n$-dimensional quasitorus
(respect. complex quasitorus) with {\em quasi-Lie algebra} $\d$ (respect.  $\dc$). 
The quasitorus $D$ (respect. $D_{\SC}$) is a real (respect. complex) 
quasifold covered by one chart.
The corresponding projection $\d\rightarrow D$ (respect.
$\dc\rightarrow\Dc$) is the {\em exponential mapping} and we denote it by 
$\exp$ (respect. $\exp_{\SC}$).}
\end{defn}
Notice that when $Q$ is a true lattice then the quasitorus $\d/Q$ is a torus.
\begin{example}\label{quasicerchio}{\rm
Consider the quasilattice $\Z+\alpha\Z$ in $\R$, with 
$\alpha\in\R\setminus\Q$.
A basic example of real quasifold is E. Prato's 
quasicircle: the real quasitorus $D^1_{\alpha}=\R/\Z+\alpha\Z$
\cite{p}. Notice that, if $\alpha$ is taken in $\Q$, 
then $D^1_{\alpha}$ 
is either an orbifold or $S^1$.
A basic example of complex quasifold is the complexification 
$\left(D^1_{\alpha}\right)_{\SC}$ of
$D^1_{\alpha}$, this is the complex quasitorus given by 
$\C/(\Z+\alpha\Z)$.}\end{example}
\subsection{Complex stratifications}\label{sezionestratificazioni}
We recall the definition of decomposition and stratification by  
quasifolds (cf. \cite{ns}).
\begin{defn}\label{decomposition}{\rm
Let $X$ be a  topological space.
A {\it decomposition of $X$ by quasifolds} is a collection of
disjoint, locally closed, connected 
quasifolds ${\cal T}_{F}$ ($F\in{\cal F}$),
called {\it pieces}, such that
\begin{enumerate}
\item The set $\cal F$ is finite, partially ordered and 
has a maximal element;
\item $X=\bigcup_F{\cal T}_F$;
\item ${\cal T}_F\cap{\overline{\cal T}}_{F'}\neq\emptyset$ iff ${\cal
T}_F\subseteq{\overline {\cal T}}_{F'}$ iff $F\leq F'$;
\item the piece corresponding to the maximal element is open and dense in $X$.
\end{enumerate}
The space $X$ is then said to be an {\em $m$-dimensional
space decomposed by quasifolds}, where $m$ is the 
dimension of the maximal piece. 
We call the maximal piece the regular piece and the other pieces  singular.}  
\end{defn}
A mapping from a decomposed space $X$ to a decomposed space $X'$ is smooth 
(respect. a diffeomorphism) if it is a continuos mapping (respect. a 
homeomorphism) that respects the decomposition and is smooth 
(respect. a diffeomorphism) when restricted to pieces.

A stratification is a decomposition that satisfies a local triviality 
condition.

Let  $L$ be a space decomposed by quasifolds.
A {\it cone over $L$}, denoted by $C(L)$, is the space $[0,1)\times L/\sim$, 
where two points $(t,l)$ and $(t',l')$ in $[0,1)\times L$
are equivalent if and only if $t=t'=0$. The decomposition of $L$
induces a decomposition of the cone.
Now let $t$ be a point in a quasifold $\cal T$ 
and let $B\simeq \tilde{B}/\Gamma$ 
a local model of $\cal T$ containing $t$. 
The decomposition of $L$ induces a decomposition of
the product $\tilde{B}\times C(L)$. 
Suppose, in addition, that $\G$ acts freely on $\tilde{B}$ and that
the space $L$ is endowed with an action of $\Gamma$ that preserves the 
decomposition; then the product $\tilde{B}\times C(L)$ is acted on by $\G$ and
the quotient $(\tilde{B}\times C(L))/\Gamma$ inherits the 
decomposition of $\tilde{B}\times C(L)$. Moreover
the quotient $(\tilde{B}\times C(L))/\Gamma$ fibers over $B$ with fiber
$C(L)$. 
\begin{defn}\label{stratificazione}{\rm Let $X$ be an $m$-dimensional
space decomposed by real quasifolds, the decomposition of  $X$
is said to be a {\it stratification by quasifolds} if each
singular piece ${\cal T}$, called {\it stratum}, satisfies the
following conditions:
\begin{enumerate}
\item[1.] let $r$ be the real dimension of $\cal T$,
for every point $t\in{\cal T}$ there exist: 
an open neighborhood $U$ of $t$ in $X$; 
a local model $B\simeq\tilde{B}/\Gamma$ in $\cal T$ containing 
$t$ and such that $\G$ acts freely on $\tilde{B}$;
a $(m-r-1)$-dimensional compact space $L$,  called the {\it link} of $t$, 
decomposed by quasifolds;
an action of the group $\Gamma$ on $L$, 
preserving the decomposition of $L$ 
and such that the pieces of the induced decomposition of
$\tilde{B}\times C(L)/\Gamma$ are  quasifolds;
finally a homeomorphism
$h\,\colon\,(\tilde{B}\times C(L))/\Gamma\lorw U$ that
respects the decompositions and  
takes each piece of $(\tilde{B}\times C(L))/\Gamma$
diffeomorphically into the corresponding piece of $U$;
\item[2.] 
the decomposition of the link $L$ satisfies condition 1.
\end{enumerate}
}\end{defn}
The definition is recursive and, since the dimension of $L$ decreases
at each step, we end up, after a finite number of steps,
with links that are compact quasifolds.
\begin{remark}\label{quandobanale}{\rm Notice that, 
if the discrete groups $\G$'s are finite 
for any possible $F$, $t\in{\cal T}_F$ and $B$, then
the twisted products  $\tilde{B}\times\C(L)/\Gamma$
become trivial and the singular strata turn out to be smooth manifolds, 
since $\G$'s act freeely. Therefore our stratification satisfies in this case 
the local triviality condition of the classical definition of 
stratification, morover, strata are smooth, with the only possile exception
of the principal stratum, that might be an orbifold.
}\end{remark}
We can then ask that the stratification be endowed with
a global complex structure \cite{cr}, namely strata are complex quasifolds and
the complex structure of each stratum is compatible with the stratification,
these are very strong requirements that are not usually satisfied, 
for example complex stratifications are usually far from being holomorphically
locally trivial, however, toric varieties and our toric spaces satisfy the 
following definition: 
\begin{defn}[Complex stratification]\label{cstratification}{\rm The 
stratified space $X$ is endowed with a complex structure
if the following conditions are satisfied:
\begin{itemize}
\item[1.] for each link $L$ there exist: a compact space $Y$ stratified
by quasifolds; a smooth surjective map $s: L\lorw Y$; a $1$-parameter
subgroup $S$ of a real torus, acting smoothly on $L$, with 
$0$-dimensional stabilizer, such that, for each $y\in Y$, the fiber 
$s^{-1}(y)$ is diffeomorphic to the quotient of $S$ by the stabilizer 
of $S$ on the fiber itself;
\item[2.]  each piece of the stratified spaces $X$, $C(L)$'s and 
$Y$'s is endowed with a complex structure;
\item[3.] the natural projection $C(L)\setminus\{\hbox{cone pt}\}\lorw Y$,
induced by $s$, when restricted to each piece, is holomorphic,
with fiber over each $y\in Y$ biholomorphic to the quotient of  
$S_{\SC}$ by the stabilizer on the fibre itself; 
\item[4.] the space $X$ is locally biholomorphic to  the product
$\tilde{B}\times C(L)/\Gamma$, that is the identification mapping $h$
is a biholomorphism when restricted to pieces.
We call $Y$ {\em complex link}  and $C(L)$ {\em complex cone} over $Y$.
\end{itemize}
}\end{defn}
\section{The construction of the quotient}\label{quotient}
In this section we describe the subsequent steps that lead to the construction
of the quotient space associated to a convex polytope $\D$.
\subsection{The open subset $\C^d_{\D}$ in $\C^d$}\label{opensubset}
Let $\d$ be a real vector space of dimension $n$, and let $\Delta$
be a convex polytope of dimension $n$ in the dual space $\ddu$.
Let $d$ be the number of facets of $\D$, write $\Delta$ as intersection of 
half spaces 
\begin{equation}\label{polydecomp}
\D=\bigcap_{j=1}^d\{\;\mu\in\ddu\;|\;\langle\mu,X_j\rangle\geq\lambda_j\;\}
\end{equation}
where $\xd$ are the chosen generators of the $1$-dimensional cones
of the fan dual to $\D$; the coefficients
$\lambda_j$'s are uniquely determined by the $X_j$'s.  
For each open face $F$ of $\D$ we denote by $I_F$ the subset of
$\{1,\ldots,d\}$ such that
\begin{equation}\label{facce}
F=\{\,\mu\in\Delta\;|\;\langle\mu,X_j\rangle=\lambda_j\;
\hbox{ if and only if}\; j\in I_F\,\}.
\end{equation}
The $n$-dimensional open face of $\D$ corresponds to the empty set.
A partial order on the set of all faces of $\D$ is defined by
setting $F\leq F'$ (we say $F$ contained in $F'$) if
$F\subseteq\overline{F'}$. The polytope $\Delta$ is the disjoint
union of its faces.
Let $r_F=\hbox{card}(I_F)$; we have the
following definitions:
\begin{defn}{\rm A $p$-dimensional face $F$ of the polytope is
said to be {\em singular} if $r_F>n-p$,
{\em regular} if $r_F=n-p$.}
\end{defn}
\begin{remark}{\rm Let $F$ be a $p$-dimensional singular face in $\d^*$, 
then $p<n-2$.
For example: a polytope in $(\R^2)^*$ is simple; the singular faces of a 
nonsimple polytope in $(\R^3)^*$ are $0$-dimensional. 
}\end{remark}
Now let $\K$ be one of the following sets $\C,\C^*,\R,\R^*$,
all of them considered naturally immersed in $\C$.
Let $J$ be a subset of $\{1,\dots,d\}$ and let $J^c$ be its complement.
We denote by
\begin{equation}\label{notazioneprincipale}
\K^J=\{\zd\in\C^d\;|\;z_j\in\K\;\;\hbox{if}\;\; j\in J,\;
z_j=0\;\;\hbox{if}\;\; j\notin J
\}.\end{equation}
We have $\K^d=\K^J\times(\K)^{J^c}$.
Let $\vz\in\K^d$, we denote by $\vz_J$ its projection onto the factor
$\K^J$. By $T^J$ we denote the subtorus $\{(t_1,\dots,t_d)\in T^d\;|\:
t_j=1\;\hbox{if}\;j\notin J\}$.
Let $F$ be a $p$-dimensional face and let $I_F$ be the corresponding
set of indices. To lighten the notation we shall omit the $I$ and simply
write $\K^F$, $\K^{F^c}$, $T^F$,...,instead of $\K^{I_F}$, 
$\K^{I_{F^c}}$, $T^{I_F}$.

Let us denote by $\C^d_{\D}$ the open subset of $\C^d$ given by
\begin{equation}\label{apertopulito}
\C^d_{\D}=\cup_{F\in\D} \C^F\times {\C^*}^{F^c}
\end{equation}
Notice that in the definition of the open subset $\C^d_{\D}$
only the combinatorics of the polytope intervenes. Moreover, the open subset
$\C^d_{\D}$ coincides with the one defined in \cite{cox} for the rational 
case. 
\subsection{The group $N_{\SC}$}
It is in the definition of the group acting on $\C^d_{\D}$ that
nonrationality comes in. In order to obtain the group we
adopt a generalization to nonrational polytopes \cite{p} of the
Delzant procedure \cite{delzant}.
Let $Q$ be a quasilattice in the space $\d$
containing the elements $X_j$ (for example 
$\hbox{Span}_{\Z}\{X_1,\cdots,X_d\}$) 
and let $\{\ed\}$ denote the standard
basis of $\rd$; consider the surjective linear mappings
\begin{equation}\label{pi}
\begin{array}{cccc}
\pi \,\colon\,& \R^d & \lorw & \d\\
    &   e_j& \longmapsto & X_j,
\end{array}
\qquad
\begin{array}{cccc}
\pic \,\colon\,& \C^d & \lorw & \dc\\
    &   e_j& \longmapsto & X_j,
\end{array}
\end{equation}
when no ambiguity can arise we shall drop the subscript $\C$.
Consider the quasitorus $\d/Q$ and its complexification
$\dc/Q$. Each of the mappings $\pi$ and $\pic$ induces a group 
homomorphism, $$\Pi
\,\colon\, T^d=\R^d/\Z^d\lorw \d/Q$$ and $$\Pic \,\colon\, 
\Tdc=\C^d/\Z^d\lorw \dc/Q.$$
We define $N$ to be the kernel of the mapping $\Pi$ and $\Nc$ to be the kernel 
of the mapping $\Pic$. 
The mapping $\Pic$ defines the isomorphism
\begin{equation}\label{qtiso}
\Tdc/\Nc\longrightarrow \dc/Q
\end{equation}
\begin{remark}\label{polar}{\rm As in the simple case \cite{cx}, we have, 
for the complexified group $\Nc$, the polar decomposition, namely:
\begin{equation}\label{polareq}
\Nc=NA,
\end{equation}
where $A=\exp(i\n)$: every element $w\in\Nc$ can be written 
uniquely as $x\,\exp(i Y)$ where $x\in N$ and $Y\in\n$. This follows from the
definition of $N$ and $\Nc$, indeed
$\Nc=\{\,\exp(Z)\;|\;Z\in\C^d\quad\hbox{and}\quad\pic(Z)\in Q\,\}$. 
Write $Z=X+iY$, then $\pic(Z)\in Q$ if and only if $\pi(X)\in Q$ and $\pi(Y)=0$, which implies
(\ref{polareq}).}
\end{remark}
Notice that neither $N$ nor $\Nc$ is a torus unless $Q$ is a lattice.

The next statement will be of great help when
analysing the local structure of our space.
Define ${\cal I}$ to be the set of subsets $I$ of $\{1,\dots,d\}$ 
such that $\{X_j\;|\;j\in I\}$ is 
a basis of  $\d$ and $I\subset I_{\mu}$ for a vertex $\mu$ of $\D$. 
If a vertex $\mu$ is non singular then ${I}_{\mu}$ itself is in $\cal I$.
Define the group $$\G_I=N\cap T^I,$$ 
observe that $\G_I$ is discrete by (\ref{exact}).
We are now able to state a Lemma, which was proved in the simple case in
\cite[Lemma 2.3]{cx}. The proof goes through with no changes, we  
briefly recall it for completeness.
\begin{lemma}  
\label{toro} Let $I\in{\cal I}$.
Then we have
that\\ 1. $T_{\SC}^d/T_{\SC}^I\simeq N_{\SC}/\Gamma_{I}$;\\ 
2. $N_{\SC}=\Gamma_{I}\exp{(\n+i\n)}$;\\ 3. given any
complement $\b_{\SC}$ of $\C^I$ in $\C^d$, we have that
$$\n_{\SC}=\{\,Z-\pi_{I}^{-1}(\pi(Z))\;|\;Z\in\b_{\SC}\,\}.$$
\end{lemma}
\proof 1. The natural group homomorphism 
$
N_{\SC}\lorw T_{\SC}^d/T_{\SC}^I
$ induces an isomorphism
$T_{\SC}^d/T_{\SC}^I\simeq N_{\SC}/\Gamma_I$.\\
2. Every element in $N_{\SC}$ can 
be written
in the form $\exp{(Z)}$, where $Z\in\C^d$ is such that 
$\pi(Z)\in Q$. Then
$Z-\pi_{I}^{-1}(\pi(Z))\in\n_{\SC}$, and 
$\exp{(\pi_{I}^{-1}(\pi(Z)))}\in\G_{I}$. The group $\G_{I}\cap\exp{(\n_{\SC})}$ is
not necessarily trivial, so the decomposition need not to be unique.
\\ 3. Split $V\in\n_{\SC}$ as
$V=Z_1+Z_2$ according to $\C^d=\C^I\oplus\b_{\SC}$, then
$\pi(V)=0$ implies that $Z_1=-\pi_{I}^{-1}(\pi(Z_2))$. \qed
The following corollary is important, it means that the action of $N$ on 
$\C^d$ is not far from being proper:
\begin{cor}\label{compattezza} Let $I\in{\cal  I}$ and 
let $\{c_m\}\in N$ be a sequence, then there 
exists a subsequence $\{c_k\}$ such that for each $k$ the element 
$c_k=\gamma_k b_k$ with 
$\gamma_k\in \Gamma_I$, $b_k\in N$ and the sequence $\{b_k\}$ is
converging to an element $b\in N$.\end{cor}
\proof By Lemma~\ref{toro} the sequence $c_m$ can be written as
$$c_m=\gamma'_m \left(\exp(-\pi_I^{-1}(\pi(Y_m)))\exp(Y_m)\right),$$ with $Y_m\in\R^{I^c}$,
therefore there exists a subsequence $Y'_k\in\R^{I^c}$ such that
$Y_k-Y'_k\in\Z^{I^c}$ and $Y'_k$ is convergent in $\R^{I^c}$.
Therefore $c_k=\gamma_k b_k$ with $\gamma_k\in\Gamma_I$ and
$b_k=\exp(-\pi_I^{-1}(\pi(Y_k)))\exp(Y_k)$ convergent in $N$ by continuity.
\qed
\subsection{The $\exp(i\n)$-orbits}
Let $\vz$ be a point in $\C^d_{\D}$, we say that the $A$-orbit $A\vz$ is closed
if it is closed in $\C^d_{\D}$. The first step towards the construction of our toric space
as quotient is given by the following theorem:
\begin{thm}[Closed orbits]\label{closedorbits}
Let $\vz\in\C^d_{\D}$. Then the $A$-orbit through $\vz$, $A\vz$, is 
closed if and only if there exists a face $F$ such that $\vz$
is in $({\C^*})^{F^c}$.
Moreover, if $A\vz$ is nonclosed, then its closure contains one and only one
closed $A$-orbit.
\end{thm}
\proof  Working on the proof of \cite[Thm~2.1]{cox} did help us to find
an argument for the proof of the first part of our statement.
Let $\vz\in ({\C^*})^{F^c}$. We want to prove that 
$A\vz$ is closed in $\C^d_{\D}$. Observe first that 
$A\vz\subset(\C^*)^{F^c}$.
Since $\C^d_{\D}=\cup_{G}(\C^G\times ({\C^*})^{G^c})$, 
it suffices to prove that 
$A\vz$
is closed in the open subset $\C^G\times ({\C^*})^{G^c}$ for each face $G$ 
of $\D$ such that 
$A\vz\subset \C^G\times ({\C^*})^{G^c}$, namely for each face $G$ such that
$I_F\subseteq I_G$, that is such that $G\subseteq F$.
Let us consider first the case of $G$ properly contained in $F$.
Let $\xi\in F\setminus G$ and $\eta\in G$: the coefficients
$$c_j=\langle \xi,X_j\rangle-\langle \eta,X_j\rangle$$
have the following properties:
\begin{itemize}
\item[] $c_j=\lambda_j-\lambda_j=0,\qquad\hbox{for}\quad j\in I_F$
\item[] $c_j=\langle \xi,X_j\rangle-\lambda_j>0
\quad\hbox{for}\quad j\in I_G\setminus I_F.$
\end{itemize}
We define on the subset $\C^{G}\times (\C^*)^{G^c}$
the continous function
\begin{equation}\label{polinomio}
P(\vw)=\Pi_{j=1}^{d}|w_j|^{c_j}
\end{equation}
The function $P$ is well defined. We prove that $P$ is invariant
under the action of $N_{\SC}$. The function $P$ is clearly invariant under
the action of $N$, we need to prove that it is invariant under the action of 
$A=\exp(i\n)$. 
Consider the exact sequence:
\begin{equation}
\label{exact}
0\lorw\n\stackrel{\iota}{\lorw}\R^d\stackrel{\pi}{\lorw}\d\lorw0
\end{equation}
and the dual sequence
\begin{equation}
\label{exactdual}
0\lorw\d^*\stackrel{\pi^*}{\lorw}(\R^d)^*\stackrel{\iota^*}{\lorw}\n^*\lorw0
\end{equation}
Let $X\in\n$ and $a=\exp(i(\iota(X))$, then
$$a\vw=(e^{-2\pi\langle \iota(X),e_1^*\rangle}w_1,\cdots,
e^{-2\pi\langle \iota(X),e_d^*\rangle}w_d)$$
and
$$P(a\vw)=e^{-2\pi\sum_{j=1}^d(\langle\xi,X_j\rangle-\langle
\eta,X_j\rangle)\langle \iota(X),e^*_j\rangle}\Pi_{k=1}^d|w_k|^{c_k}$$
Notice that
$$\sum_{j=1}^d(\langle\xi,X_j\rangle-\langle
\eta,X_j\rangle)\langle \iota(X),e^*_j\rangle=
\langle X,\iota^*(\sum_{j=1}^{d}\langle \pi^*(\xi-\eta),e_j\rangle e_j^*)
\rangle.$$
Therefore, since 
$\sum_{j=1}^{d}\langle \pi^*(\xi-\eta),e_j\rangle e_j^*\in\Im\,\pi^*$,
by (\ref{exactdual}) we obtain
$$i^*\left(\sum_{j=1}^{d}\langle \pi^*(\xi-\eta),e_j\rangle e_j^*\right)=0$$
and therefore
$$P(a\vw)=P(\vw).$$
Remark that $P(\vz)\neq0$ and define the subset
$P_{\vz}=\{\vw\in\C^{G}\times (\C^*)^{G^c}\,|\,P(\vw)=P(\vz)\}$. 
The set $P_{\vz}$ is
closed, is contained in $\C^{F}\times (\C^*)^{F^c}$ and the orbit
$A\vz$ is contained in $P_{\vz}$, since the function
$P$ is $A$--invariant.
It follows that $\overline{A\vz}$ is contained in $P_{\vz}$ and therefore in 
$\C^{F}\times (\C^*)^{F^c}$. We prove that this implies 
$$\overline{A\vz}=A\vz.$$ 
Let $\vw$ be a point in the closure of $A\vz$ in $\C^{G}\times (\C^*)^{G^c}$, 
then $\vw\in\C^{F}\times (\C^*)^{F^c}$ and
there exists a sequence $Y_n\in\n$ such that the sequence
$\exp(i Y_n)\vz$ converges to $\vw$. Therefore, for $j\notin I_F$,
$$\lim_{n\rightarrow+\infty}e^{-2\pi\langle \iota(Y_n),e^*_j\rangle}z_j=w_j.$$ 
This implies that there exists $Y'\in\R^{F^c}$ such that
$$\lim_{n\rightarrow+\infty}e^{-2\pi\langle \iota(Y_n),e^*_j\rangle}=
e^{-2\pi\langle Y',e^*_j\rangle}$$ and 
$$w_j=e^{-2\pi\langle Y',e^*_j\rangle}z_j.$$
Now remark that 
$(\pi\circ\iota)(Y_n)=0$, which gives
\begin{equation}\label{sopra}
\sum_{j\in I_F}\langle \iota(Y_n),e^*_j\rangle X_j=
-\sum_{j\notin I_F}\langle \iota(Y_n),e^*_j\rangle X_j.
\end{equation}
Let $$\d_F=\hbox{Span}\{X_j\;|\;j\in F\}.$$
The sequence 
$\sum_{j\in I_F}\langle \iota(Y_n),e^*_j\rangle X_j$ is in $\d_F$,
by (\ref{sopra}) it converges to
$-\sum_{j\notin I_F}\langle Y',e^*_j\rangle X_j$, which must
therefore lie in $\d_F$. Let $Y''\in\R^{F}$ such that
$$\sum_{j\in I_F}\langle Y'',e^*_j\rangle X_j=
-\sum_{j\notin I_F}\langle Y',e^*_j\rangle X_j,$$ and set
$Y=Y'+Y''\in\R^d$. It follows that $\pi(Y)=0$, therefore $Y\in\n$ and
$$\vw=\exp{iY}\vz.$$  
The above argument also proves that $A\vz$ is closed in 
$\C^F\times(\C^*)^{F^c}$.

Now consider a point $\vz$ such that $\vz\notin (\C^*)^{G^c}$ for
any face $G$ of $\D$. Here suggestions for the proof 
come from an argument by I. Musson \cite{musson} used by Cox in his proof. 
Take a face $F$ such that $$\vz\in \C^F\times(\C^*)^{F^c}.$$ 
Then $F$ must
be singular, otherwise we would have $\vz$ contained in 
$(\C^*)^{G^c}$ for some face $G$ containing $F$.
Therefore there exist coefficients $y_k$ not all zero such that
$$\sum_{k\in I_F}y_k X_k=0.$$
We want to prove that the coefficients $y_k$ can be chosen so that 
there exists a $j\in I_F$ with the property that:
$$
z_j\neq0,\quad y_j>0
$$
Let $I_{\vz}=\{k\in I_F\;|\;z_k=0\}$
and take 
\begin{equation}\label{facciagiusta}
E=\cap_{k\in I_{\vz}}\{\xi\in\D\;|\;\langle\xi,X_k\rangle=\lambda_k\},
\end{equation}
then either $E=F$ or $F\subset E$: if $E=F$ then $\{X_k\;|\;k\in I_{\vz}\}$
span $\d_F$, therefore we can find an index $j$ with the 
required properties; if $F\subset E$ then any 
$j\in I_E\setminus I_{\vz}$, which is nonempty by our hypothesis on $\vz$,
satisfies the required properties.
Let $y_k$ be the chosen coefficients for $k\in I_F$, set
$y_k=0$ for $k\notin I_F$ and let $$Y=(y_1,\dots,y_d),$$ 
since $\pi(Y)=0$ we have that 
$Y\in\n$. Remark now that 
$$\lim_{t\rightarrow+\infty}\exp(itY)\vz\notin A\vz,$$
the orbit is therefore nonclosed. 

Now we want to prove that $\overline{A\vz}$ contains one and only one
closed orbit. Let $E$ be the face defined by (\ref{facciagiusta}),
then $I_{\vz}\subset I_E$, moreover $\{X_k\;|\;k\in I_{\vz}\}$ span
$\d_E$. For each $j\in I_E\setminus I_{\vz}$ we can therefore find
$y_k^j\in\R$ such that:
$$\begin{array}{l}
y_j^j=1\\
y^j_k=0,\quad k\notin I_{\vz},\;\;k\neq j\\
\sum_{k=1}^d y_k^j X_k=0
\end{array}
$$
Let $Y^j=(y^j_1,\dots,y^j_d)$, we have that $\pi(Y^j)=0$, therefore $Y^j\in\n$.
Therefore
$$\lim_{t\rightarrow+\infty}
\left(\Pi_{j\in I_E\setminus I_{\vz}}\exp(itY^j)\right)\vz=\vz_{E^c}$$
implies that the closure of the orbit $A\vz$ contains the orbit 
$A\vz_{E^c}$, which is closed.

Suppose now that there is another closed orbit, $A\vu$, in $\overline{A\vz}$.
Then we have $\vu\in \C^{G^c}$ for some face $G$ such that $F\subseteq G$.
If $E\neq G$ we can take $\xi\in E\setminus E\cap G$ and
$\eta\in G\setminus E\cap G$. We
 can then construct a function $P$ as in
(\ref{polinomio}) with coefficients
$$c_j=\langle \xi,X_j\rangle-\langle \eta,X_j\rangle.$$
We have
\begin{itemize}
\item[] $c_j=\lambda_j-\lambda_j=0,\qquad\hbox{for}\quad j\in I_E\cap I_G$
\item[] $c_j=\langle \xi,X_j\rangle-\lambda_j>0
\quad\hbox{for}\quad j\in I_G\setminus (I_E\cap I_G).$
\item[] $c_j=\lambda_j-\langle \eta,X_j\rangle<0
\quad\hbox{for}\quad j\in I_E\setminus (I_E\cap I_G).$
\end{itemize}
Remark that $I_{\vz}\subset I_E\cap I_G$.
The function $P$  is well defined on the orbits $A\vz$ and $A\vu$,
it is $A$-invariant and 
therefore constant and nonzero on the $A$-orbit $A\vz$. This implies that
$I_G=I_G\cap I_E$. Taking $P^{-1}$ proves that $I_E=I_G\cap I_E$. It follows
that $E=G$.

In order to procede we rely on the following statement that will be proved
later on, in Corollary~\ref{unique}.
\begin{itemize}
\item[]
{\em Let $G$ be a face of $\D$ and let $\vw\in (\C^*)^G$, then there exist
a unique point $\vw^0$ in $A\vw$ and a unique point $\xi\in G$ such that
$$\langle \xi,X_j\rangle-\lambda_j=|w^0_j|^2,\quad j=1,\dots,d$$}
\end{itemize}
Now let $A\vw$ and $A\vu$ be two closed orbits in  $\overline{A\vz}$ and let
$(\vw^0,\xi)$ and $(\vu^0,\eta)$ as in the statement right above. Suppose
that the two orbits do not coincide, therefore $(\vw,\xi)\neq(\vu,\eta)$.
Let $P$ be the function given in (\ref{polinomio}), with coefficients
$c_j=\langle \xi,X_j\rangle-\langle \eta,X_j\rangle.$
The function $P$ is well is constant and nonzero on 
$A\vz$ and hence by continuity on the orbits $A\vu$ and $A\vw$. It follows that
$P(\vu^0)/P(\vw^0)=1$. Now let $J^{+}=\{j\;|\;c_j>0\}$ and
$J^-=\{j\;|\;c_j<0\}$, we have
$$1=P(\vu^0)/P(\vw^0)=(\Pi_{k\in J^+}|u^0_k/w^0_k|^{c_k})
(\Pi_{k\in J^-}|u^0_k/w^0_k|^{c_k})$$
now notice that
$$
c_j=|w^0_j|^2-|u^0_j|^2
$$
which leads to a contradiction, therefore $J^+=J^-=\emptyset$
and $|w^0_j|^2=|u^0_j|^2$ for
all $j\in\{1,\dots,d\}$, which implies $\vw^0=\vu^0$, since they lie
in the closure of the same $A-$orbit.
\qed
\subsection{The quotient by the nonclosed group $N_{\SC}$}
Theorem~\ref{closedorbits} allows us to define on the open set $\C^d_{\D}$ the
following equivalence relation: two points $\vz$ and $\vw$ are equivalent
with respect to the action of the group $N_{\SC}$,
\begin{equation}\label{equivalenzasegno}
\vz\sim_{\! N}\vw,
\end{equation}
if and only if 
\begin{equation}\label{equivalenza}
\left(N(\overline{A\vz})\right)\cap\left(\overline{A\vw}
\right)\neq\emptyset,
\end{equation}
where the closure is meant in $\C^d_{\D}$.
\begin{prop}\label{propequivalenza} The relation defined by (\ref{equivalenza}) is an equivalence relation.
\end{prop}
\proof 
The relation defined in (\ref{equivalenza}) is symmetric:
let $\vz$ and $\vw$ such that $\vz\sim_{\! N}\vw$, then 
$\left(N(\overline{A\vz})\right)\cap\left(\overline{A\vw}\right)\neq\emptyset$.
Remark first that $\overline{A\vz}$ is a union of $A-$orbits and that 
$N(\overline{A\vz})$ is also a union of $A-$orbits: for the first claim take
$\vu\in\overline{A\vz}$, then there exists a sequence $a_k\in A$ such that
$a_k\vz\lorw\vu$, then, since the action is continuos, for each $a\in A$
$a_ka\vz\lorw a\vu$, therefore $A\vu\subset\overline{A\vz}$; for the second 
claim take $\vu\in N(\overline{A\vz})$ then there exists $c\in N$ and 
$\vu'\in\overline{A\vz}$ such that $\vu=c\vu'$, but, for what we have just 
observed, the whole orbit $A\vu'$ is contained in $\overline{A\vz}$, therefore
the whole orbit $A\vu$ is contained in $N(\overline{A\vz})$.
It follows that the intersection 
$\left(N(\overline{A\vz})\right)\cap\left(\overline{A\vw}\right)$ 
is a union of $A-$orbits. Let 
$A\vu\subset\left(N(\overline{A\vz})\right)\cap\left(\overline{A\vw}\right)$, 
then
there exists $c\in N$ such that $A\vu$ is an $A$-orbit in 
$c(\overline{A\vz})$, 
therefore $A(c^{-1}\vu)$ is an $A-$orbit in 
$\overline{A\vz}\cap c(\overline{A\vw})$, hence
$\left(N(\overline{A\vw})\right)\cap\left(\overline{A\vz}\right)\neq\emptyset$.
This proves that the relation is symmetric. 
We prove now that the relation is transitive.
Remark first that if
$A\vu$ is in the intersection 
$\left(N(\overline{A\vw})\right)\cap\left(\overline{A\vz}\right)\neq\emptyset$,
then its closure also lies in the intersection,
by Theorem~\ref{closedorbits}
$\overline{A\vu}$ contains one and only one closed $A$-orbit.
Let $\vz\sim_{\!N}\vv$ and 
$\vv\sim_{\! N}\vw$. The closed orbit in
$N(\overline{A\vz})\cap\overline{A\vv}$ and the 
closed orbit
in  $N(\overline{A\vw})\cap\overline{A\vv}$  
coincide since they both lie
in $\overline{A\vv}$, therefore $\vz\sim_{\! N}\vw$.
\qed
We define the space $X_{\D}$ to be the quotient of $\C^d_{\D}$ by the 
equivalence relation just defined, we denote the quotient by
$$X_{\D}=\C^d_{\D}//N_{\SC}.$$
Notice that, if the polytope is simple, then 
$\C^d_{\D}=\cup_{F\in\D}(\C^*)^F$ and 
the quotient $X_{\D}$ is just the orbit space endowed with the quotient 
topology.
\begin{remark}\label{ancherazionale}{\rm The construction of the quotient
can be carried out for any convex polytope, rational or not. If we start
with a rational polytope in a lattice $L$ and we choose the primitive 
generators of the $1$-dimensional cone in the dual fan, then the above quotient
is the toric variety associated to the polytope in $L$ \cite{cox}. In general
there are many toric spaces associated to a given convex polytope, depending on
the choice of quasilattice and generators. See the model example of the unit 
interval and the family of toric spaces associated to it \cite{p}. 
There are applications in which it is natural to consider a rational polytope 
in a quasilattice \cite{rhombus,3d}.}  
\end{remark}
\section{The structure of the quotient $X_{\D}$}
\label{decomposizionedix}
\subsection{The decomposition and the structure of the pieces}
The decomposition in pieces of the quotient $X_{\D}$ reflects the 
geometry of the polytope $\D$. 
The indexing set $\cal F$ for the decomposition of $X_{\D}$ is given by 
the set of the singular faces of the 
polytope $\D$, with the partial order defined in Section~\ref{opensubset},
with the addition of a maximal element. 
The maximal piece is
$${\cal T}_{\hbox{\scriptsize{max}}}=\cup_{F\,\hbox{\scriptsize{reg}}}
(\C^F\times(\C^*)^{F^c})/N_{\SC}=
\cup_{F\,\hbox{\scriptsize{reg}}}({\C^*})^{F^c}/N_{\SC}.$$
Then there is a piece ${\cal T}_F$ for each singular face $F$ of $\D$:
$${\cal T}_F=
\{z\in\C^d_{\D}\;|\;\overline{Az}\cap(\C^*)^{F^c}\neq\emptyset\}//N_{\SC}.$$
Theorem~\ref{closedorbits} implies that
the space $X_{\D}$ is given by the union of the regular and singular pieces.
The structure of $X_{\D}$ as decomposed space and the properties that 
characterize $X_{\D}$ as a toric space associated to $\D$ are described in 
this section: 
\begin{thm}[Quasifold structure of strata]
\label{stratiquasifold}
The subset ${\cal T}_F$ of $X_{\D}$ corresponding to a $p$-dimensional
singular face of $\D$ is a $p$-dimensional complex quasifold.
The subset ${\cal T}_{\hbox{\scriptsize{max}}}$
is an $n$-dimensional complex quasifold. These subsets give a
decomposition by complex quasifolds of $X_{\D}$.
\end{thm}
\proof {\em Regular piece}
Consider the regular piece ${\cal T}_{\hbox{\scriptsize{max}}}$. 
The proof that
it is an $n$-dimensional complex quasifold goes very similarly
to the proof, given in \cite[Thm2.2]{cx}, that the space correponding
to a simple polytope, with a given choice of normals and quasilattice,
is an $n$-dimensional complex quasifold.  
Notice that $\cup_{F\,\hbox{\scriptsize{reg}}}({\C^*})^{F^c}$ is covered by the
sets  $\widehat{V}_I=(\C^I\times(\C^*)^{I^c})\cap\left(\cup_{F\,
\hbox{\scriptsize{reg}}}({\C^*})^{F^c}\right)$,
with $I\in{\cal I}$. 
Let $\tilde{V_I}\subset\C^I$ be the image of the natural projection 
mapping from $\left(\C^I\times(\C^*)^{I^c}\right)
\cap\left(\cup_{F\;
\hbox{\scriptsize{reg}}}({\C^*})^{F^c}\right)$
to $\C^I$.
Consider the mapping
$$
p_I:\tilde{V_I}\lorw\widehat{V}_I$$ defined by
$$\left(p_I(\vz)\right)_j=\left\{\begin{array}{l}
z_j\;\hbox{if}\;j\in I\\
1\;\hbox{if}\;j\in I^c
\end{array}\right.$$
The mapping $p_I$ induces a homeomorphism
$$\phi_I: \tilde{V}_I/\G_I\lorw\widehat{V}_I/N_{\SC}.$$
The proof that $\phi_I$ is a homeomorphism, as well as the proof that the 
charts \linebreak $(\widehat{V}_I/N_{\SC},\phi_I,\tilde{V}_I/\G_I)$  give an atlas
of the regular piece, go along the lines of the proof of \cite[Thm.2.2]{cx}, 
we will not repeat the argument here.

{\em Singular pieces}
Let $F$ be a singular face,
we want to characterize those points $\vz$ which are equivalent, with respect
to (\ref{equivalenza}), to points in $(\C^*)^{F^c}$.
Let $I_h$, with $h=1,\dots,m_F$, be the subsets of $I_F$ such that
$$F=\cap_{j\in I_h}\{\xi\in\D\;\:\;\langle \xi,X_j\rangle=\lambda_j\}$$
From the proof of Theorem~\ref{closedorbits} it follows that 
$$\{\vz\in\C^d_{\D}\;|\;\overline{A\vz}\cap(\C^*)^{F^c}\neq\emptyset\}=
\cup_{h=1}^{m_F}(\C^*)^{I_h^c}.$$
Moreover, if we set $\pi_h:(\C^*)^{I_h^c}\lorw(\C^*)^{F^c}$, for a given
$N_{\SC}$-invariant subset $U$ of $(\C^*)^{F^c}$ we have that:
\begin{equation}\label{utile}
\{\vz\in\C^d_{\D}\;|\;\overline{A\vz}\cap U\neq\emptyset\}=
\cup_{h=1}^{m_F}\pi_h^{-1}(U).\end{equation}
This, together with Theorem~\ref{closedorbits}, implies that
the inclusion mapping
$$(\C^*)^{F^c}\hookrightarrow\cup_{h=1}^{m_F}(\C^*)^{I_h^c}$$
induces a homeomorphism when passing to the quotient by $N_{\SC}$, where
on the left we consider the geometric quotient and on the right the quotient
by the equivalence relation (\ref{equivalenza}). 
We can therefore identify the piece ${\cal T}_F$ with the orbit space
$(\C^*)^{F^c}/N_{\SC}$.
To conclude the proof of the theorem we need to show that the 
$(\C^*)^{F^c}/N_{\SC}$ is a complex quasifold of dimension $p$,
where $p$ is the dimension of $F$.
Choose an $I\in{\cal I}$ such that $\hbox{card}(I_F\cap I)=n-p$ and
denote by
\begin{equation}\label{gammacheck}
\check{\Gamma}_I=\G_I/(\G_I\cap T^F).
\end{equation}
The action of $\check{\Gamma}_I$ on $(\C^*)^{I\setminus I\cap I_F}$ is
well defined. 
Now consider the mapping
$p_F:(\C^*)^{I\setminus I\cap I_F}\lorw (\C^*)^{F^c}$
defined by
$$\left(p_F(\vz)\right)_j=\left\{\begin{array}{l}
z_j\;\hbox{if}\;j\in I\setminus I\cap I_F\\
1\;\hbox{if}\;j\in (I\cup I_F)^c
\\0\;\hbox{if}\;j\in I_F\end{array}\right.$$
The mapping $p_F$ induces a bijective continuos mapping 
\begin{equation}\label{fieffe}
\phi_F:(\C^*)^{I\setminus I\cap I_F}/{\check \G}_{I}
\lorw (\C^*)^{F^c}/N_{\SC}.\end{equation}
The mapping $\phi_F$ is also open since
the mapping
$$\begin{array}{ccc}\C^{{(I\cup I_F)^c}}\times (\C^*)^{I\setminus I\cap I_F}&
\lorw& (\C^*)^{F^c}\\
(\vw,\vz)&\longmapsto&
(\exp(\vw)\exp(\pi_I^{-1}(\pi(\vw)))p_F(\vz)\end{array}$$
is not only surjective but has surjective differential at every point.
Therefore $\phi_F$ is a homeomorphism and 
$(\C^*)^{I\setminus I\cap I_F}/{\check \G}_{I}$ gives a chart covering 
the $p$-dimensional complex quasifold ${\cal T}_F$.
\qed
\begin{cor}\label{stratisingolari}
The singular stratum ${\cal T}_F$ corresponding to a singular face $F$
can be identified with the orbit space $(\C^*)^{F^c}/N_{\SC}$, which
is precisely, by (\ref{qtiso}), the $D_{\SC}$-orbit corresponding to the face 
$F$. The maximal stratum is the union of the orbits of
$D_{\SC}$ corresponding to the regular faces, in particular, it contains
the orbit corresponding to the interior of the polytope.
\end{cor}
\proof See proof of Theorem~\ref{stratiquasifold}.
\begin{cor}\label{mappaaperta} The projection mapping 
$\C^d_{\D}\lorw \C^d_{\D}//N_{\SC}$ is open.\end{cor}
\proof The statement can be easily proved by making use of
(\ref{utile}).
\begin{prop}\label{denseorbit} 
The $n$-dimensional complex quasitorus $D_{\SC}$ acts continuosly
on $X$, with a dense open orbit. Moreover the
restriction of the $D_{\SC}$-action
to each piece of the space $X$ is
holomorphic.
\end{prop}
\proof The proof is a simple consequence of Corollary~\ref{stratisingolari},
see also the proof of the analogous result in \cite{cx,ns}. The dense open
orbit is the one corresponding to the interior of the polytope.
\subsection{Building blocks: links and complex cones}
Consider the singular $p$-dimensional piece ${\cal T}_F$. We show that all of the points of
${\cal T}_F$ have the same link. We first describe this link on the polytope $\D$.
Let
$\jmath_F:\d_F\hookrightarrow\d$ be the 
inclusion mapping.
Define
$$\Sigma^{\diamond}_F=\bigcap_{j\in I_F}\{\;\mu\in\ddu\;|\;
\langle\mu,X_j\rangle\geq\lambda_j\;\}.$$ 
By 
projecting $\Sigma^{\diamond}_F$ onto $\d^*_F$ we obtain the cone
\begin{equation}\label{sigmaeffe}\Sigma_F=\jmath_F^*(\Sigma^{\diamond}_F),\end{equation} 
it is an 
$(n-p)$-dimensional cone in $\d^*_F$ with vertex $\jmath_F^*(F)$.
If $G$ is a $q$-dimensional face of $\D$ containing $F$, 
then $\jmath_F^*(G)$ is a $(q-p)$-dimensional face of ${\Sigma_F}$.
By slicing the cone $\Sigma_F$ with a hyperplane transversal to its 
faces we obtain a polytope $\D_F$ which is the link of the face $F$. Each face
of $\D$ containing $F$ gives rise to a face of $\D_F$.
More precisely: for each $j\in I_F$ we can find an $s_j\in(0,1]$ such that,
taken $X_0=\sum_{j\in I_F}s_j X_j$, the intersection
$$\D_F={\Sigma_F}\cap\{\xi\in\d_F^*\;|\;\langle\xi,X_0\rangle=
\sum_{j\in I_F}\lambda_j s _j+1\}$$ is a nonempty convex 
polytope of dimension $(n-p-1)$.  
Let $G$ be a $q$-dimensional face of $\D$ 
properly containing $F$, then 
$$\jmath_F^*(G)\cap\{\xi\in\d_F^*\;|\;\langle\xi,X_0\rangle=
\sum_{j\in I_F}\lambda_j s_j+1\}$$
is a $(q-p-1)$-dimensional face of $\D_{F}$, which 
is singular in $\D_{F}$ if and only if $G$ is singular in $\D$.
The complex spaces corresponding to these newly defined
convex polyhedra are the building blocks of our stratified space.
For each point in ${\cal T}_F$ the cone $C(L)$ of 
Definition~\ref{stratificazione} 
will be the complex space corresponding to the 
polyhedral cone  $\Sigma_F$, whilst the space $Y$ of  
Definition~\ref{cstratification} will be
the complex space corresponding to $\D_F$.
Let us now construct these spaces. Recall that
$\d_F=\hbox{Span}\{X_j\;|\;j\in F\}$.
Notice first that 
$Q_F=\d_F\cap Q$ is a quasilattice.
Consider the convex polyhedral cone $\Sigma_F$, 
together with the normals $X_j$, with $j\in I_F$, and the quasilattice
$Q_F$.
We have the short exact sequence
\begin{equation}\label{esattapersigmaf}
\xymatrix{
0\ar[r]&\n^F\ar[r]^{\iota_F}&\R^F\ar[r]^{\pi_F}&\d_F\ar[r]&0}
\end{equation}
then
$$N^F=\Ker(T^{F}\lorw\d_F/Q_F).$$
Remark that
$N^F=N\cap T^F$ with $\hbox{dim}\,N^F=r_F-n+p$.
With the procedure described in Section~3 we construct the quotient
$$C(L_F)=
\C^{F}_{\Sigma_F}//N^F_{\SC}=\C^{F}//N^F_{\SC}$$
this is a space decomposed by complex quasifolds,
the decomposition is the one
induced by that of $\Sigma_F$ in regular and singular faces. Remark that
$\Sigma_F$ is not a polytope, it is a polyhedral cone, but the construction 
described in the previous sections applies with no changes; notice also
$\C^F_{\Sigma_F}=\C^F$.

Now let $\hbox{ann}(X_0)$ be the annihilator of $X_0$ and let 
$$\xymatrix{
\hbox{ann}(X_0)\ar@{^{(}->}[r]^-{k_F}&(\d_F)^*}$$ 
be the natural inclusion.
Then $k^*_F$ projects $\d_F$ onto
the $(n-p-1)$-dimensional space $(\hbox{ann}(X_0))^*\simeq 
\d_F/\langle X_0\rangle$.
We continue to denote by $\D_F$ the polytope 
$\D_{F}$ viewed in the subspace 
$\hbox{ann}(X_0)$:
\begin{equation}\label{deltaeffe}\D_{F}=
\bigcap_{j\in I_F}\{\;\xi\in\hbox{ann}(X_0)\;|\;
\langle\xi,k_F^*(X_j)\rangle\geq\lambda_j-\langle\xi_0,X_j\rangle\;\}
\end{equation} where $\xi_0$ is a point in the affine hyperplane with which
we cut $\Sigma_F$. 
Now  apply
the construction described in Section~3 to
the polytope $\D_{F}$, 
with the choice of normals $k_F^*(X_j)$ and quasilattice 
$Q_{F,0}=k_F^*(\d_F\cap Q)$.
We obtain the exact sequence
\begin{equation}\label{esattaperdeltaf}
\xymatrix{
0\ar[r]&\n^F_0\ar[r]^{\iota_{F,0}}&\R^F\ar[r]^<<<<<{k_F^*\circ\pi_F}
&(\hbox{ann}(X_0))^*\ar[r]&0}
\end{equation}
Let
$\s=\hbox{Span}\{s'_1,\dots,s'_d\}$ where
$s'_j=s_j$ if $j\in I_F$ and $s'_j=0$ otherwise.
It is easy
to check that
$$\n^F_0=\n^F\oplus\s.$$
The subgroup of $T^{F}$ that we need is then
$$N_0^F=\Ker\left(T^{F}\lorw(\hbox{ann}(X_0))^*/Q_{F,0}\right).$$
Remark that
$$N_0^F/N^F\simeq\hbox{exp}(\s).$$
Notice that the polytope $\D_F$, obtained by
slicing the cone $\Sigma_F$, is combinatorially equivalent to
$\Sigma_F\setminus\{\hbox{cone point}\}$,
therefore $\C^{F}_{\D_F}$ does not depend on the choice of  $X_0$, whilst
the group $N_0^F$ does.
The space corresponding to $\D_F$ with this set of data is the space 
$X_{\D_F}$, decomposed by complex quasifolds, given by the quotient
$$\C^{F}_{\D_F}//(N_0^F)_{\SC}$$
\begin{lemma} The natural mapping
$$s_F:\C^{F}_{\Sigma_F}//N^F_{\SC}\setminus\{[0]\}
\lorw\C^{F}_{\D_F}//(N_0^F)_{\SC}$$
satisfy Definition~\ref{cstratification}.\end{lemma}
\proof We have already observed that 
$(N_0^F)_{\SC}/N_{\SC}^F\simeq
\exp(\s_{\SC})$,
therefore
the continuos mapping $s_F$ is surjective, with
fibre isomorphic to  the quotient of the complex group
$\exp(\s_{\SC})$ by the stabilizer of $\exp(\s_{\SC})$ on the fiber
itself. 
Moreover $s$ naturally respects the decomposition and is holomorphic when
restricted to each piece.\qed
\begin{remark}{\rm The coefficients $s_j$ can always be chosen
in such a way that the group $\exp(\s)$ is compact or acts freely 
on $\C^F$, nonetheless, since from the examples 
it is clear that there are choices of $\s$ that are
natural (see \cite[Example~3.6]{ns}),
we prefer to have freedom of choice on the coefficients and
to leave Definition~\ref{cstratification} as it is,
with no stronger requirements on the action of $S$.
}\end{remark}
\subsection{Holomorphic local triviality}
Let us first define and analyse closely the candidate product space.
The action of the finitely generated group
$\check{\G}_I$ 
on the product 
$(\C^*)^{I\setminus(I\cap I_F)}
\times(\C^{F}//N^F_{\SC})$
is defined by
$[\gamma](\vz,[\vw])=(\gamma\vz,[\gamma\vw])$.
Observe that the action is well defined and it is free on the first factor.
Now define on the product $(\C^*)^{I\setminus(I\cap I_F)}\times\C^{F}$
the following equivalence relation:
\begin{equation}\label{sim1}
(\vw_1,\vz_1)\sim_1(\vw_2,\vz_2)
\end{equation}
if and only if there exists $\gamma\in\G_I$ such that 
$$\gamma\vw_1=\vw_2$$
and
$$\gamma\vz_1\sim_{\!N^F}\vz_2.$$
It is easy to check that 
\begin{equation}\label{prima}
\left((\C^*)^{I\setminus(I\cap I_F)}
\times\C^{F}//N^F_{\SC}\right)/
\check{\G}_{I}\simeq(\C^*)^{I\setminus(I\cap I_F)}
\times\C^{F}/\sim_1
\end{equation}
Consider now the action of the quotient group $N_{\SC}/N^F_{\SC}$
on the product
$(\C^*)^{F^c}
\times(\C^{F}//N^F_{\SC})$
defined by 
$$[g]\cdot(\vw,[\vz])=(g\vw,[g\vz])$$
where $g\in N_{\SC}$, $\vw\in(\C^*)^{F^c}$ and $\vz\in\C^F$.
On $(\C^*)^{F^c}\times\C^F$ define the following equivalence relation:
\begin{equation}\label{sim2}
(\vw_1,\vz_1)\sim_2(\vw_2,\vz_2)
\end{equation}
if and only if there exists $g\in N_{\SC}$ such that
$$g\vw_1=\vw_2$$
and
$$g\vz_1\sim_{\!N^F}\vz_2.$$
It is easy to check that
\begin{equation}\label{seconda}
\left((\C^*)^{F^c}
\times\C^{F}//N^F_{\SC}\right)/(N_{\SC}/N^F_{\SC})\simeq
\left((\C^*)^{F^c}
\times\C^{F}\right)/\sim_2
\end{equation}
We want to prove the following 
\begin{lemma}\label{tildetilde} The spaces
$$
\left((\C^*)^{I\setminus(I\cap I_F)}
\times\C^{F}\right)/\sim_1$$ and 
$$
\left((\C^*)^{F^c}
\times\C^{F}\right)/\sim_2
$$ 
are diffeomorhic as decomposed spaces, moreover the diffemorphism is a 
biholomorphism when restricted to the pieces.
\end{lemma}
\proof First of all remark that, since the projection 
$\C^F\lorw \C^F//N^F_{\SC}$ is open by Corollary~\ref{mappaaperta},
we have that the projections
$$(\C^*)^{I\setminus(I\cap I_F)}
\times\C^{F}\lorw\left((\C^*)^{I\setminus(I\cap I_F)}
\times\C^{F}\right)/\sim_1$$ and 
$$
(\C^*)^{F^c}
\times\C^{F}\lorw
\left((\C^*)^{F^c}
\times\C^{F}\right)/\sim_2$$
are both open.
The candidate diffeomorphism between the spaces in the statement 
is the mapping $f$ defined by
$f([\vw,[\vz]])=[\vw+\underline{1},[\vz]]$ where
$\vw\in\C^{I\setminus I\cap I_F}$, $\vz\in\C^F$ and 
$\underline{1}\in \C^{(I\cap I_F)^c}$ is defined by
\begin{equation}\label{uno}
\underline{1}_j=1\quad\hbox{for}\; 
j\notin I\cup I_F\qquad
\underline{1}_j=0\quad\hbox{for}\; 
j\in I\cup I_F.\end{equation}
It is easy to check that $f$ is bijective, in particular injectivity is
due to the fact that $$N_{\SC}\cap T^{I\cup I_F}=\G_I N^F_{\SC}.$$
It is also a straightforward check to prove that $f$ is continuous.
The key point is to prove that $f$ is open. 
Let
$U$ be an open subset of $(\C^*)^{F^c}\times\C^F$, saturated with
respect to $\sim_2$ and let $W\times\underline{1}\times V$ be contained
in $U$, with $W$ an open subset 
of 
$(\C^*)^{I\setminus I\cap I_F}$ and $V$ an open subset of 
$\C^F$,  invariant under the
actions of $N^F$ and $\G_I$.
In order to prove that $f$ is open it suffices to prove that
$$N_{\SC}\left(W\times\underline{1}\times V\right)$$ is open (it is obviously
contained in $U$).
Consider the mapping
\begin{equation}\label{aperta}
\begin{array}{ccc}
W\times V\times\C^{(I\cap I_F)^c}&\lorw&(\C^*)^{F^c}\times\C^F\\
((\vw,\vz),V)&\longmapsto&\exp(V)\exp(\pi_I^{-1}(\pi(V))(\vw,\underline{1},\vz)
\end{array}\end{equation}
the image of this mapping is a subset of 
$N_{\SC}(W\times\underline{1}\times V)$, since
by Lemma~\ref{toro} $$\exp(V)\exp(\pi_I^{-1}(\pi(V))\in N_{\SC}.$$ We prove
that it is exactly $N_{\SC}(W\times\underline{1}\times V )$: take an element 
$g\in N_{\SC}$, then, by Lemma~\ref{toro}, there exist $\gamma\in\G_I$ and
$Z\in\C^{I^c}$ such that $g=\gamma\exp(Z)\exp(\pi^{-1}_I(\pi(Z))$. 
Split $Z=Z_1+Z_2$ according to the decomposition 
$\C^{I^c}=\C^{(I\cup I_F)^c}\times\C^{I_F\setminus(I_F\cap I)}$, then
$$g=\gamma\exp(Z_1)\exp(\pi_I^{-1}(\pi(Z_1)))
\exp(Z_2)\exp(\pi_I^{-1}(\pi(Z_2)).$$
Now observe that $$\exp(Z_2)\exp(\pi_I^{-1}(\pi(Z_2))=
\exp(Z_2)\exp(\pi_{I\cap I_F}^{-1}(\pi_F(Z_2))$$ is an element of $N^F_{\SC}$
and recall that the open subset $V$ is invariant by the action of 
$N_{\SC}^F$ and by the action of $\G_I$. Therefore
$g(\vw,\underline{1},\vz)=
\exp(Z_1)\exp(\pi_I^{-1}(\pi(Z_1))(\vw',\underline{1},\vz')$ with
$(\vw',\vz')\in W\times V$. This proves our assertion that the saturated of
$W\times\underline{1}\times V$ under the action of $N_{\SC}$ is exactly 
the image of the mapping
(\ref{aperta}), on the other hand it is easy to check that 
the differential of the 
mapping (\ref{aperta}) is never zero, its image is therefore open.
It is now straightforward to check that the restriction of $f$ to the pieces
is a biholomorphism.\qed

We are now ready to prove the local holomorphic triviality of our decomosition:
\begin{lemma}[Holomorphic local triviality] Recall that 
${\cal T}_F\simeq (\C^*)^{I\setminus I\cap I_F}/{\check \G}_{I}$
and that  $((\C^*)^{F^c}\times\C^F)//N_{\SC}$ is an open subset of
$X_{\D}$ containing
${\cal T}_F$.
There exists a mapping $h_F$ from
$\left((\C^*)^{I\setminus(I\cap I_F)}
\times\C^{F}//N^F_{\SC}\right)/
\check{\G}_{I}$ onto $((\C^*)^{F^c}\times\C^F)//N_{\SC}$,
which is a homeomorphism and a biholomorphism restricted to the pieces of
the respective decompositions. 
\end{lemma}
\proof  
By (\ref{prima}, \ref{seconda}) and Lemma~\ref{tildetilde} 
it suffices to prove that the mapping
$$\begin{array}{ccccc}h_F&:&
(\C^*)^{F^c}\times\C^F
/\sim_2
&\lorw&((\C^*)^{F^c}\times\C^F)//N_{\SC}\end{array}$$
is a homeomorphism and a biholomorphism restricted to the pieces of
the respective decompositions. 
Consider two points
$(\vw,\vz)$ and $(\vw',\vz')$ in $(\C^*)^{F^c}\times\C^F$. 
Then $(\vw',\vz')\sim_N(\vw,\vz)$ if and only
if there exists an element $g\in N_{\SC}$ and a face $G$ in $\D$ such that
$(g(\vw'+\vz'))_{G^c}=(\vw,\vz)_{G^c}$, if and only if 
$g\vw'=\vw$ and $(g\vz')_{G^c}=\vz_{G^c}$ if and only if
$g\vw'=\vw$ and $g\vz'\sim_{N^F}\vz$ if and only if 
$(\vw',\vz')\sim_2(\vw,\vz)$. It is then easy to check that
the mapping $h_F$ is a homeomorphism and a
biholomorphism when restricted to pieces.\qed
\subsection{The interplay with the symplectic setup}
In this section we will prove the statement that was used in the 
proof of Theorem~\ref{closedorbits}.
Moreover, in order to conclude that $\cquo$ is a complex stratified space, 
we need to prove
that the spaces $\C^F//N^F_{\SC}$ are cones over a compact stratified space
$L$ satisfying the definition of link, moreover all of the mappings 
involved have to satisfy the requirements of the 
Definitions~\ref{stratificazione},\ref{cstratification}.
General results, for example on bundles, are not of immediate application 
to our spaces because of their topology; therefore,
although a description of the building bloks of our stratification can 
be given within our set up (see the first row of diagram~\ref{diagramma3}), 
we make use of the interplay with the symplectic quotients in order to give a 
neat description
of the link $L_F$ and of $X_{\Sigma_F}$ as a real cone over it. Hence,
before going on to describe the cone $X_{{\Sigma}_F}$, we will 
briefly recall from \cite{ns} the symplectic construction for the nonsimple
case (cf. \cite{p} for the simple case). 

Consider the convex sets $\D$, $\Sigma_F$ and $\D_F$.
We have already defined the groups $N$, $N^F$ and $N^F_0=N^F\exp(\s)$.
We briefly recall the construction of the moment mappings relative to these
polyhedral sets. 
A moment mapping with respect to  
the standard action of the torus $T^d$ on $\C^d$ is given by
\begin{equation}\label{upsilon}
\Upsilon(\vz)=\sum_{j=1}^d
(|z_j|^2+\lambda_j)e_j^*,
\end{equation} 
where the $\lambda_j$'s are given in
(\ref{polydecomp}).
Analogously we consider the mappings on $\C^F$ given by
\begin{equation}\label{upsiloneffe}
\Upsilon_F(\vz)=\sum_{j\in I_F}
(|z_j|^2+\lambda_j)e_j^*,
\end{equation} 
and
\begin{equation}\label{upsiloneffezero}
\Upsilon_{F,0}(\vz)=\sum_{j\in I_F}
(|z_j|^2+\lambda'_j)e_j^*,
\end{equation} 
with $\lambda_j'=\lambda_j-\langle \xi_0,X_j\rangle$.

A moment mapping with respect to the induced Hamiltonian action of  $N$ 
on $\C^d$ is then: 
$\Psi_{\D}\,\colon\,\C^d\rightarrow \n^*$ given by 
$$\Psi_{\D}={\iota}^*\circ\Upsilon.$$
In the same manner the mappings
$$\Psi_{\Sigma_F}=\iota_F^*\circ\Upsilon_F$$
and
$$\Psi_{\D_F}=\iota_{F,0}^*\circ\Upsilon_{F,0}$$
are moment mappings relative to the
Hamiltonian action of the $(r_{F}-n+p)$-dimensional group $N^{F}$
on ${\C}^{F}$ 
and to the
Hamiltonian action of the $(r_{F}-n+p+1)$-dimensional group $N^{F}_0$
on ${\C}^{F}$ respectively.
The explicit expression of these moment mappings
are
$$\begin{array}{ccccc}
\Psi_{\D}&:&\C^d&\lorw&\n^*\\
&&\vz&\longmapsto&\sum_{j=1}^d\left(|z_j|^2+\lambda_j\right)\iota^*(e^*_j)
\end{array}
$$
\begin{equation}\label{esplicito2}
\begin{array}{ccccc}
\Psi_{\Sigma_F}&:&\C^F&\lorw&\n_F^*\\
&&\vz&\longmapsto&\sum_{j\in I_F}|z_j|^2\iota^*_F(e_j^*)
\end{array}
\end{equation}
and finally
\begin{equation}
\label{esplicito3}
\begin{array}{ccccc}
\Psi_{\D_F}&:&\C^F&\lorw&\n_F^*\oplus\R\\
&&\vz&\longmapsto&(\sum_{j\in I_F}|z_j|^2\iota^*_F(e_j^*),
\sum_{j\in I_F}s_j|z_j|^2-1)
\end{array}
\end{equation}
Where in (\ref{esplicito3}) we have identified 
$(\n\oplus\s)^*\simeq\n^*\oplus\R$.
The symplectic quotient corresponding to $\D$ 
is then the reduced space 
$$M_{\D}=\Psi^{-1}_{\D}(0)/N,$$
whilst for $\Sigma_F$ and $\D_F$ we have
$$M_{\Sigma_F}=\Psi_{\Sigma_F}^{-1}(0)/N^F$$ and 
$$M_{\D_F}=\Psi_{\D_F}^{-1}(0)/N^F_0$$ respectively. 
These symplectic quotients
are spaces stratified by symplectic quasifolds, they are
endowed with the continuos effective Hamiltonian action of
the quasitori $\d/Q$, $\d_F/Q_F$ and $(\hbox{ann}(X_0))^*/Q_{F,0}$ 
respectively. 
The quasitori act smoothly on the strata, with
moment mappings given by the restriction to the strata of the mappings 
$\Phi_{\D}=(\pi^*)^{-1}\circ\Upsilon$, $\Phi_{\Sigma_F}=(\pi_F^*)^{-1}\circ
\Upsilon_F$ and
$\Phi_{\D_F}=(\pi_F^*\circ k_F)^{-1}\circ\Upsilon_{F,0}$ respectively.
The constants in (\ref{upsilon},\ref{upsiloneffe},\ref{upsiloneffezero}) 
have been chosen in such a way that
the images of the mappings $\Phi_{\D}$, $\Phi_{\Sigma_F}$
and $\Phi_{\D_F}$ are respectively $\D$, $\Sigma_F$ and $\D_F$ \cite{ns}.

In order to relate the complex and symplectic quotients we need 
to prove the following statement, which is an 
adaptation to our case of the results contained in \cite[Appendix 1]{g}, to 
which we refer for further details.
\begin{lemma}\label{opencone}
The mapping $\Psi_{\D}:\C^d\lorw\n^*$ satisfies the following properties:
\begin{itemize}
\item[1.]  $\Psi_{\D}^{-1}(0)\subset \cup_{F}(\C^*)^{F^c}$; 
morover, for each face $F$, the set
$(\C^*)^{F^c}$ intersects $\Psi_{\D}^{-1}(0)$ in at least one point;
\item[2.] let $\vz$ be any point in $\C^d_{\D}$, then $\Psi(A\vz)$ is a cone, 
which is open as a subset of a suitable linear subspace 
of $\n^*$; the cone depends only on the set $I_{\vz}=\{j\;|\;z_j=0\}$;
\item[3.] let $\vz\in \C^d_{\D}$, then there exists an 
$a\in A$ such that $\Psi_{\D}(a\vz)=0$ if and only if $\vz\in(\C^*)^{F^c}$
for some face $F$ of $\D$;
\item[4.] $\Psi_{\D}$ is a proper mapping, in particular $\Psi_{\D}^{-1}(0)$ 
is compact.
\end{itemize}
\end{lemma}
\proof
Let $\vz\in\C^d$. We have that $\Psi(\vz)=0$ if and only if
$\iota^*(\sum_{j=1}^d(|z_j|^2+\lambda_j)e_j^*)=0$ if and only if, 
by (\ref{exactdual}), there
exists a (unique) point $\xi\in\d^*$ such that 
$$\pi^*(\xi)=\sum_{j=1}^d(|z_j|^2+\lambda_j)e_j^*$$
if and only if for every $k=1,\cdots,d$ we have 
$$\langle \xi,X_k\rangle=|z_k|^2+\lambda_k,$$
if and only if $\xi\in\D$. This implies point 1.

Denote by $\alpha_k$ the elements of $\n^*$ given by 
$\alpha_k=-2\pi\iota^*(e_k^*)$. Consider the stabilizer of $A$ at $\vz$,
given by $\exp(i(\n\cap\R^{I_{\vz}}))$ and let $\r$ be the orthogonal 
complement
of $\n\cap\R^{I_{\vz}}$ in $\n$.
Finally consider the subset of $\n^*$ given by $\hbox{Span}\{\alpha_k\;|\;
k\notin I_{\vz}\}$. Now remark that
$$\hbox{Span}\{\alpha_k\;|\;k\notin I_{\vz}\}\simeq
\hbox{ann}(\n\cap\R^{I_{\vz}})\simeq\r^*$$
where $\hbox{ann}(\n\cap\R^{I_{\vz}})$ denotes the annihilator of
$\n\cap\R^{I_{\vz}}$.
The subset $\Psi(A\vz)$
of $\n^*$ can be identified with the image of the mapping
\begin{equation}
\label{image}
\begin{array}{ccc}
\r&\lorw&\r^*\\
Y&\longmapsto&
\sum_{k\notin I_{\vz}}e^{\alpha_k(Y)}|z_k|^2\alpha_k+
\lambda\end{array}
\end{equation}
where $\lambda=\sum_{k=1}^{d}\lambda_k\iota^*(e^*_k)$.
Remark now that the mapping above is the Legendre
transform of the function
$$\begin{array}{ccccc}F_{\vz}&:&\r&\lorw&\R\\
&&Y&\longmapsto&\sum_{k\notin I_{\vz}}e^{\alpha_k(Y)}|z_k|^2+\lambda(Y)\end{array}.$$
It is easy to check that the Hessian of $F_{\vz}$ is positive definite, since 
the $\alpha_k$'s generate $\r^*$, therefore $F_{\vz}$ is strictly 
convex. This implies that the image of the mapping (\ref{image}), namely
$\Psi(A\vz)$, is the open
convex cone in $\hbox{Span}\{\alpha_k\;|\;k\notin I_{\vz}\}$ given by
$$\{\sum_{k\notin I_{\vz}}t_k\alpha_k+\lambda\;|\;t_k>0\}$$
moreover $\r$ is mapped diffeomorphically onto the open cone by the 
mapping (\ref{image}). This proves point 2.

Now we are able to prove point 3: let $\vz\in(\C^*)^{F^c}$ for some $F$ in
$\D$. By point
(i) there exists a point $\vw$ in $(\C^*)^{F^c}$ such that the $A$-orbit
$A\vw$ intersects the zero set $\Psi^{-1}(0)$, therefore the cone in
$\n^*$ corresponding to the orbit $\vw$ must contain $0$. But by point (ii)
this holds also for the orbit $A\vz$. By the same kind of argument we can 
deduce that, if $\vz$ is not in $(\C^*)^{F^c}$ for some face $F$, then
the $A$-orbit $A\vz$ does not intersect the zero set $\Psi^{-1}(0)$.

In order to prove the last point of our statement, let us first prove that
the zero set is compact. As we have already remarked $\Psi(\vz)=0$ if
and only if $\Upsilon(\vz)\in \hbox{Ker}(\iota^*)$. Therefore
$\Psi^{-1}(0)=\Upsilon^{-1}(\Im(\Upsilon)\cap \hbox{Ker}\iota^*)$.
But $\Im(\Upsilon)\cap \hbox{Ker}\iota^*=
\{\xi\in(\R^d)^*\;|\;\langle\xi,e_k\rangle\geq\lambda_k,\;\;k=1\cdots,d\}\cap\hbox{Im}(\pi^*)$.
This can be deduced from the definition of $\Upsilon$ together with
(\ref{exactdual}). Therefore  $\Im(\Upsilon)\cap \hbox{Ker}(\iota^*)$ is exactly
$\pi^*(\D)$, in particular it is compact. This implies, since $\Upsilon$ is
proper, that $\Psi^{-1}(0)$ is compact. 

Observe now that, since $\iota^*$ is a linear projection and 
$\hbox{Im}(\Upsilon)$ is the positive orthant shifted by $\lambda$, we have
that $\hbox{Im}(\Upsilon)\cap (\iota^*)^{-1}(\eta)$
is compact for 
every $\eta\in\n^*$. And again properness of $\Upsilon$ implies 
that $\Psi$ is also proper. 
\qed
\begin{cor}\label{unique} 
Let $F$ be a face of $\D$, let $\r$ be the orthogonal complement
of $\n\cap\R^F$ in $\n$ and let $\vz$ be a point in $(\C^*)^{F^c}$. 
Then there exists a unique point $\vx\in\zset$, a unique point
$\xi\in F$ and a unique $Y\in\r$ such that
$$A\vz\cap\zset=\{\vx\},$$
$$\exp(iY)\vz=\vx$$
and
$$\langle \xi,X_j\rangle-\lambda_j=|x_j|^2,\quad j=1,\dots,d$$
\end{cor}
\proof The argument is the same used in \cite[Remark~3.3]{cx}.
Suppose we have two points of intersection: $\vx_1$ and $\vx_2$, then
there is an $a\in A$ such that $a\vx_1=\vx_2$, but either $a$ is in the 
stabilizer of $\vx_1$, and therefore $\vx_1=\vx_2$, or $a$ is not in the 
stabilizer of $\vx_1$. In this case it moves $\vx_1$ out of $\zset$, 
contradiction.\qed 
\begin{cor}\label{unique2} 
Let $\vz$ be a point in $\C^d_{\D}$, then there exists a unique 
point $\vx\in\zset$ such that $\overline{A\vz}\cap\zset=\{\vx\}$.\end{cor}
\proof The statement follows by Corollary~\ref{unique} and by 
Theorem~\ref{closedorbits}
\qed 
The inclusion mapping $\Psi_{\D}^{-1}(0)\hookrightarrow\C^d_{\D}$ induces
a mapping 
\begin{equation}\label{lachi}
\chi_{\D}\,\colon\,\Psi_{\D}^{-1}(0)/N\lorw\cquo
\end{equation}
between the 
two quotients. By Corollary~\ref{unique2} we can also define the 
surjective mapping
\begin{equation}\label{llatildechi}
\begin{array}{ccccc}
\Xi_{\Delta}&:&\C^d_{\D}&\lorw&\Psi^{-1}(0)\\
&&\vz&\lorw&\overline{A\vz}\cap\Psi^{-1}(0).
\end{array}
\end{equation}
This implies that $\chi_{\D}$ is bijective and 
$\chi_{\D}^{-1}([\vz])=[\Xi_{\D}(\vz)]$.
\begin{lemma}\label{lift} Suppose that the mapping $\chi_{\D}$ is a 
homeomorphism and
take a converging  sequence $\vz_n\lorw\vz$ in $\C^d_{\D}$, then the sequence
$\Xi_{\D}(\vz_n)$ converges to $\Xi_{\D}(\vz)$.
\end{lemma} 
The proof is direct and is left to the reader.
\begin{lemma}\label{iso1}
Let ${\cal T}_F$ be a singular stratum, corresponding to the 
singular face $F$. The complex and symplectic quotients that express the
singular stratum in the symplectic and complex setup are naturally isomorphic.
\end{lemma}
\proof Let $I$ be such that card$(I\cap I_F)=n-p$ and
let $a_{ij}$ be the $d\times n$ matrix of the projection 
$\pi:\R^d\lorw \d$ with respect to the standard basis and the basis
$\{X_j\;|\;j\in I\}$.
Then recall that the stratum ${\cal T}_F$, given by the quotient $\left(\Psi^{-1}(0)\cap (\C^*)^{F^c}\right)/N$,
can be identified with the quotient
$\tilde{B}_F/\check{\G}_{I\setminus (I\cap I_F)}$, where
$$\tilde{B}_F=\{\vz\in(\C^*)^{I\setminus (I\cap I_F)}\;|\;
\sum_{h\in I}a_{hk}(|z_{h}|^2+\lambda_{h})-\lambda_k)>0,\;k\in (I\cup I_F)^c\}.$$

We have to prove that the mapping
$(\chi_F)_{\scriptstyle{loc}}$ in
 the following diagram is a diffeomorphism:
\begin{equation}
\xymatrix{
(\C^*)^{I\setminus (I\cap I_F)}/\check{\Gamma}_{I\setminus (I\cap I_F)}
\ar[r]^-{\phi_F}&(\C^*)^{F^c}/N_{\SC}\\
B/\check{\Gamma}_{I\setminus (I\cap I_F)}\ar[u]^{(\chi_F)_{\scriptstyle{loc}}}
\ar[r]^-{\phi_F^s}&\Psi^{-1}(0)\cap(\C^*)^{F^c}/N\ar[u]^{\chi_F}
}
\end{equation}
The mapping $\chi_F$ is induced
by the inclusion of $\Psi^{-1}(0)\cap(\C^*)^{F^c}$ into $(\C^*)^{F^c}$;
the homeomorphism $\phi_F$ was defined in (\ref{fieffe}); the homeomorphism
$\phi_F^s$ is defined by
$\phi_F^s([x])=[x+\tilde{x}]$, where $\tilde{x}\in \C^{(I\cup I_F)^c}$ is 
given by
$$\tilde{x}_k=\sqrt{
\sum_{h\in I}a_{hk}(|x_{h}|^2+\lambda_{h})-\lambda_k)}$$
\cite{ns};
The mapping 
$(\chi_F)_{\scriptstyle{loc}}$ is the diffeomorphism
that sends $\vz\in B$ to 
$$-\exp(\pi_I^{-1}(\pi(X))_{I\setminus(I\cap I_F)}\vx$$
where $$X_k=-\frac{1}{2\pi}\,log(\sum_{h\in I}a_{hk}(|x_{h}|^2+\lambda_{h})-\lambda_k)$$
for $k\in (I\cup I_F)^c$ and $X_k=0$ otherwise.\qed

We prove now for the symplectic case an analogous of Lemma~\ref{tildetilde}:
\begin{lemma}
The quotients
$$(\tilde{B}_F\times \Psi^{-1}(0)/N^F)/\check{\G}_{I}$$ and 
$$\left((\Psi^{-1}(0)\cap(\C^*)^{F^c})\times
\Psi_{\Sigma_F}^{-1}(0)/N^{F}\right)/(N/N_F)
$$
are diffeomorphic.
\end{lemma}
\proof The action of $N/N^F$ is defined as follows: consider
$[g]\in N/N^F$ and $(\vx,[\vz])\in (\Psi^{-1}(0)\cap(\C^*)^{F^c})\times
\Psi_{\Sigma_F}^{-1}(0)/N^{F}$, then
$[g](\vx,[\vz])=(g\vx,[g\vz])$. It easy to check that this action is
well defined. Consider now the mapping: 
$$\begin{array}{ccc}
(\tilde{B}_F\times \Psi^{-1}(0)/N^F)/\G_{I\setminus (I\cap I_F)}&\lorw& 
\Psi^{-1}(0)\cap(\C^*)^{F^c}\times
\Psi_{\Sigma_F}^{-1}(0)/N^{F}/(N/N_F)\\
\,[\vx,[\vz]]&\longmapsto&[\vx+\phi_F(\vx),[\vz]]
\end{array}
$$
It is straightforward to check that this mapping is bijective and continuos. 
The argument for the proof that it is closed is also straightforward, it is
is a direct consequence of Corollary~\ref{compattezza}.
\qed
\subsection{The main results}
We will complete in this section the proof of the two theorems below,
they will be proved together as proofs are entangled.
We shall procede by induction on the depth of the polytope $\D$, 
combining the complex and symplectic setup. As final result we obtain 
not only the proof that our decomposed space is indeed a stratification but 
also the prove that $X_{\D}$ is naturally
isomorphic to its symplectic counterpart $M_{\D}$. 
\begin{thm} [$X_{\D}$ is a stratified space]
\label{ultimo}  Let $\d$ be a real vector space of dimension $n$, and let
$\D\subset\ddu$ be a convex polytope. Choose inward-pointing normals to the
facets of $\D$, $X_1,\ldots,X_d\in\d$, and let $Q$ be a quasilattice containing
them.  The quotient $X=\C^d_{\D}//N_{\SC}$ is a complex stratified 
space, in particular
for each singular face $F$ of the polytope $\D$ there exist a 
compact link $L_F$ and a complex link $Y_F$ satisfying the requirements of 
the definitions~\ref{stratificazione},\ref{cstratification}.\end{thm}
\begin{thm} [The complex and symplectic quotients are isomorphic]
\label{teoremadellachi}
Let \linebreak $\d$ be a vector space of dimension $n$, and let
$\D\subset\ddu$ be a convex polytope. Choose inward-pointing normals to the
facets of $\D$, $X_1,\ldots,X_d\in\d$, and let $Q$ be a quasilattice containing
them. 
Then the mapping $$\chi_{\D}\,\colon\,\Psi_{\D}^{-1}(0)/N\lorw\cquo$$ is an equivariant
homeomorphism with respect to the actions of $D$ and $D_{\SC}$ respectively.
The restriction of $\chi_{\D}$
to each stratum is a diffeomorphism of quasifolds. 
Moreover the induced symplectic form on each stratum
is compatible with its complex structure, so that strata have the structure
of
 K\"ahler quasifolds.
\end{thm}
\proof \textbf{of Theorems~\ref{ultimo},\ref{teoremadellachi}}
We start by giving:
{\em
\begin{itemize}
\item a description of the real link $L_F$;
\item a characterization of the
quotient $\C^{F}//(N^{F}_{\SC})$ as a real cone over $L_F$;
\item 
a diffeomorphism between complex and symplectic quotients at the cone level.
\end{itemize}
}
We first define the notion of depth of $\D$.
Let $F$ be a singular face of the polytope $\D$. 
We call singularity depth of $F$
the minimum integer $m$ such that there exists a sequence of faces
$F_0<\cdots<F_j<\cdots<F_m$ such that $F_0=F$, $F_m$ is regular and $F_j$
is singular for all $0\leq j<m$. By definition the singularity 
depth of regular faces is set to be $0$. We define the
singularity depth of a polytope to be the maximum singularity depth
attained by its faces.
We then procede by induction on the polytope depth. Suppose that $\D$ has depth
$0$, namely $\D$ is a simple polytope. In this case, treated in \cite{cx}, 
the complex quotient is a geometric quotient and the mapping
$\chi_{\D}:\squo\lorw\C^d_{\D}/N_{\SC}$ is proved to be an equivariant 
diffeomorphism such that the induced symplectic structure on $\C^d_{\D}/N_{\SC}$ is K\"ahler. Continuity of $\chi_{\D}$ is straightforward; bijectivity
of $\chi_{\D}$ follows from Lemma~\ref{opencone}. The mapping $\chi_{\D}$ 
is then computed explicitly on the local charts and turns out to be a 
local diffeomorphism, which implies that its inverse is also continous.

Suppose now that Theorems~\ref{ultimo},\ref{teoremadellachi} hold for  
polytopes of depth less than or equal to $n-1$, we want to prove that they 
hold for polytopes of depth $n$.
Consider a polytope $\D$ of depth $n$. Let $F$ be a $p$-dimensional 
singular face and let $\Sigma_{F}$ and $\D_{F}$ be the convex sets
associated to $F$ as described in (\ref{sigmaeffe}) and (\ref{deltaeffe}).  
Notice that, in order to do so, we have to choose a suitable vector $X_0\in\d_{F}$. Consider the following diagram:
\begin{equation}\label{diagramma3}
\xymatrix{
{\scriptstyle
(\C^{F}//N^{F}_{\SC})\setminus\{
 [0]
\}}
\ar[r]^{q_2}
\ar@/^2pc/[rr]|s
\ar @{} [dr]|*+[o][F-]{2}
&
{\scriptstyle
\C^{F}_{\D_{F}}//N^{F}_{\SC}\exp(i\s)}
\ar @{} [dr]|*+[o][F-]{1}
\ar [r]^{q_1}
&
{\scriptstyle \C^{F}_{\D_{F}}//{N^{F}_0}_{\SC}}
\\
{\scriptstyle
(\Psi_{\Sigma_F}^{-1}(0)/N^{F})\setminus \{[0]\}}
 \ar[u]^{\chi_{\Sigma_F}}\ar[r]^{p_2}
\ar@/_2pc/[rr]|{s'}
&
\scriptstyle{
(\Psi_{\D_F})^{-1}(0)/N^{F}}
 \ar[u]^{\chi'_{\D_F}}\ar[r]^{p_1}&
\scriptstyle{
(\Psi_{\D_F})^{-1}(0)/N^{F}_0}
 \ar[u]^{\chi_{\D_F}}
}
\end{equation}
Consider first the diagram $\xymatrix
{*+[o][F-]{1}}
$.
The mapping $\chi_{\D_F}$ is a diffeomorphims 
by the induction hypothesis. Proposition~\ref{propequivalenza}
can be applied in order to define the quotient
$\C^{F}_{\D_{F}}//N^{F}_{\SC}\exp(i\s)$:
two points $\vz$ and $\vw$ in $\C^{F}_{\D_{F}}$ are equivalent if and only
if 
$$N^F\left(\overline{A_F\exp(i\s)\vz}\right)\cap\overline{(A_F\exp(i\s)\vw)}
\neq\emptyset.$$ It was proved  
in \cite{ns} that the quotient
$(\Psi_{\D_F})^{-1}(0)/N^{F}$ is the link and in particular
it is a space stratified by 
quasifolds of real dimension $2n-2p+1$.
The inclusion $(\Psi_{\D_F})^{-1}(0)\hookrightarrow \C^F_{\D_F}$ induces
the continuos mapping $\chi'_{\D_F}$, which is bijective 
by Lemma~\ref{opencone}, the inverse mapping being induced by the surjective 
mapping $\Xi_{\D_F}:\C^{F}_{\D_F}\lorw(\Psi_{\D_F})^{-1}(0)$.
Lemma~\ref{lift} implies that the mapping $\chi'_{\D_F}$ is a diffeomorphism.
It is easy to check that diagram $\xymatrix
{*+[o][F-]{1}}
$ is commutative. 

Consider now Diagram $\xymatrix
{*+[o][F-]{2}}
$. 
The mapping 
$q_2$ is the natural projection, the mappings $\chi_{\Sigma_F}$ 
is defined by (\ref{lachi}) 
and 
$\chi'_{\D_F}$ is defined right above.
We need to prove that $\chi_{\Sigma_F}$
is a diffeomorphism of stratified spaces.
Let $\vx\in\Psi_{\Sigma_F}^{-1}(0)$, $\vx\neq0$, we prove first that the mapping
$p_2$ that makes the diagram commute is given by:
$$p_2([\vx])=[\vx/|\vx|_{s}]$$
where
$$|\vx|_s=\sqrt{\sum_{j=1}^{d}s_j |x|^2_j}.$$
Let $\vz\in\C^{F}$ such that 
$$\{\vx\}=\Psi^{-1}_{\Sigma_F}(0)\cap\,\overline{A\vz}$$ 
therefore
$$\chi_{\Sigma_F}^{-1}([\vz])=[\vx].$$
Notice that
$\vz\in\C^F_{\D_F}$, since the set of points in $\C^F$ which are 
not equivalent to $0$ under $\sim_{N^F}$ are given exactly by $\C^F_{\D_F}$.
Recall that the closure of the orbit $\overline{A\exp(i\s)\vz}$, in
$\C^F_{\D_F}$, is a union of orbits.
From Lemma~\ref{opencone} the 
set $\Psi_{\D_F}(\overline{A\exp(i\s)\vz})$ 
in  $\n_{F}^*\oplus\R$ is given by 
the union of the cones images of these orbits.
The closure of  $\Psi_{\D_F}(\overline{A\exp(i\s)\vz})$ 
is a polyhedral cone
in  $\n_{F}^*\oplus\R$, whose
cone point 
is $(0,-1)$ by (\ref{esplicito3}). Notice that
$(0,-1)\notin\Psi_{\D_F}(\overline{A\exp(i\s)\vz})$ but 
the half line $[0,-(1-t))$, $t>0$, connecting the cone point with the 
origin, is contained in $\Psi_{\D_F}(\overline{A\exp(i\s)\vz})$.
Now observe that the points of the set $\overline{A\exp(i\s)\vz}$ 
corresponding to the half line $[0,-(1-t))$ are exactly those
given by $\overline{A\exp(i\s)\vz}\cap\Psi_{\Sigma_F}^{-1}(0)$. This implies
that $\chi^{-1}_{\Sigma_F}$ projects the orbit $A\exp(i\s)\vz$ 
onto the curve in
$\Psi^{-1}_{\Sigma_F}(0)/N^{F}$ given by $[\sqrt{t}x/|\vx|_s]$, with $t>0$.
Therefore 
$p_2([\vx])=[\vx/|\vx|_{s}]$
and
the quotient $X_{\Sigma_F}=\C^F//N^F_{\SC}$ is a real cone over 
the quotient $\C^F_{\D_F}//N^F_{\C}\exp(i\s)$, which is the link
of each point lying in the orbit ${\cal T}_F$.
Since $\Sigma_F$, without the cone point, has the same depth of $\D_F$, 
the mapping $\chi_{\Sigma_F}$ is a diffeomorphism away from the cone 
point by the induction hypothesis. The argument above shows that
the mapping $\chi_{\Sigma_F}$ extends continuously to the cone point, thus
giving a diffeomorphism of decomposed spaces 
(for a different proof see Remark~\ref{mappacontinua}).

We have thus proved all the three items above. Now we need to prove that
{\em
\begin{itemize}
\item  The complex and symplectic quotient are diffeomorphic.
\end{itemize}}
We work on the next diagram in order to check that the natural 
mapping $$\chi_{\D}:\squo\lorw\cquo,$$ 
for a polytope of depth $n$, is a 
diffeomorphism such that the symplectic structure of each stratum of  $\squo$ 
induces,  via $\chi_{\D}$, a K\"ahler structure on the corresponding stratum 
of $\cquo$.

First of all it is easy to check that 
the mapping $\chi_{\D}$ is continuos and bijective, the inverse 
mapping
$\chi^{-1}_{\D}$ is induced by $\Xi_{\D}$.
We need to prove that the inverse of $\chi_{\D}$ is also continuos. This is 
true around regular points by the quoted result on simple polytopes
\cite[Thm.]{cx}.
Consider the  following diagram, which describes the local trivializations
of the complex and symplectic quotients and their relationship:
\begin{equation}\label{diagramma4}
\xymatrix{
\scriptstyle{
\frac{(\SC^*)^{F^c}\times (\SC^F//N^F_{\SC})}
{N_{\SC}/N^F_{\SC}}
}
\ar[r]&
\scriptstyle{
\left((\SC^*)^{F^c}\times(\SC)^{F}\right)//N_{\SC}}
\\
\scriptstyle{
\frac{
\Psi^{-1}(0)\cap(\SC^*)^{F^c}\times
\Psi_{\Sigma_F}^{-1}(0)/N^{F}}{N/N_F}}\ar[u]\ar[r]^-{h'}
&
\scriptstyle{
\left(\Psi^{-1}(0)\cap(\SC^F\times\SC^*)^{F^c}\right)/N}\ar[u]^{\chi_{\D}}
}
\end{equation}
We require the diagram to be commutative, this defines 
the mapping $h'$: more precisely, given 
$(\vx,\vz)\in \left(\Psi^{-1}(0)\cap(\C^*)^{F^c}\right)\times
\Psi_{\Sigma_F}^{-1}(0)$, there exists a unique $a\in A$ such that
$a(\vx+\vz)\in\Psi^{-1}(0)$, the mapping $h'$ takes $[\vx,[\vz]]$ to
$[a(\vx+\vz)]$. In order to prove that $\chi_{\D}$ is closed it is enough
to prove that $h'$ is continuos. 

Consider a closed subset $C$ of the quotient
$(\Psi^{-1}(0)\cap(\C^F\times\C^*)^{F^c})/N$ and let
$\tilde{C}$ be the $N$-invariant closed subset of
$\Psi^{-1}(0)\cap(\C^F\times(\C^*)^{F^c})$ that projects onto $C$. 
We want to prove that the
inverse image $C_1$, via the mapping $h'$, of $C$, is closed in
$\left(
\Psi^{-1}(0)\cap(\C^*)^{F^c}\times
\Psi_{\Sigma_F}^{-1}(0)/N^{F}\right)/(N/N_F)$.
Let $\tilde{C}_1$ be
the $N$-invariant
subset of
$(\Psi^{-1}(0)\cap(\C^*)^{F^c})\times
\Psi_{\Sigma_F}^{-1}(0)$ that projects onto $C_1$.
We prove that
$\tilde{C}_1$
is closed.
Let $(\vx,\vz)$ be in the closure of $\tilde{C}_1$, then
there exists a sequence
$(\vx_n,\vz_n)\in \tilde{C}_1$ converging to $(\vx,\vz)$.
Let $G$ and $G_n$ be the faces of $\D$ such that
$\vx+\vz\in\C^{G^c}$ and $\vx_n+\vz_n\in\C^{G_n^c}$ and let $\r$ and $\r_n$
the orthogonal complement in $\n$ of $\n\cap\R^G$ and $\n\cap\R^{G_n}$ 
respectively. Then there exists a sequence $Y_n\in\r_n$ such that
$\exp(iY_n)(\vz_n+\vx_n)\in \tilde{C}$. The sequence $\exp(iY_n)(\vz_n+\vx_n)$
is in $\Psi^{-1}(0)$, therefore it admits a converging subsequence, that
we denote again by $\exp(iY_n)(\vz_n+\vx_n)$. 
Let $\vw_1+\vw_2\in\C^F\times(\C^*)^{F^c}$
be its limit and let $Y\in\r$ such that
$\exp(iY)(\vx+\vz)\in\Psi^{-1}(0)$. Notice first that for $j\notin I_F$ the sequence
$\exp(\langle Y_n,\iota^*(e_j^*)\rangle)$ converges to
$\exp(\langle Y,\iota^*(e_j^*)\rangle)$ and
$(\vw_1+\vw_2)_j=\exp(\langle Y,\iota^*(e_j^*)\rangle)((\vx+\vz)_j)$.
Observe now that there exists a face $H$ such that 
$(\vw_1+\vw_2)\in(\C^*)^{H^c}$.
We have
$$I_H\subset I_G\subset I_F.$$
Suppose that $I_H$ is a proper subset of $I_G$. Then 
take $\xi\in H\setminus G$ and $\eta\in G$: the coefficients
$$c_j=\langle \xi,X_j\rangle-\langle \eta,X_j\rangle$$
have the following properties:
\begin{itemize}
\item[a)] $c_j=\lambda_j-\lambda_j=0,\qquad\hbox{for}\quad j\in I_H$
\item[b)] $c_j=\langle \xi,X_j\rangle-\lambda_j>0
\quad\hbox{for}\quad j\in I_G\setminus I_H.$
\end{itemize}
We define on the subset $\C^{H}\times (\C^*)^{H^c}$
the continous function
$$
P(z)=\Pi_{j=1}^{d}|z_j|^{c_j}
$$
We have that
$$\lim_n\,P(\exp(iY_n)(\vx_n+\vz_n))=P(\vw_1+\vw_2)\neq0$$
on the other hand the function $P$ is invariant under the action on
$N_{\SC}$, therefore
$$\lim_n\,P(\exp(iY_n)(\vx_n+\vz_n))=
\lim_n\,P(\vx_n+\vz_n)=P(\vx+\vz)=0.$$
It follows that $H=G$ and $\vw_1+\vw_2=\exp(iY)(\vx+\vz)$. This proves
that $\vx+\vz\in \tilde{C}_1$ is closed and therefore $h'$ is continous.
\qed
\begin{cor} The space $X_{\D}$ is compact.
\end{cor}
\begin{remark}{\rm 
Theorem~\ref{teoremadellachi} improves the result of local
triviality that was found in \cite{ns}:
diagram (\ref{diagramma4}) implies that
the link {\em does not} depend on the point but just on the stratum, 
moreover the local trivialization is of the form 
$$(\tilde{B}_F\times C(L_F))/\check{\G}.$$ where $\tilde{B}_F/\check{\G}$ 
is a chart that covers the whole stratum ${\cal T}_F$}.\end{remark} 
\begin{remark}\label{mappacontinua}{\rm We can date back to
the work \cite{kempfness} the study of the
relationship between symplectic and complex quotients. 
We have given a direct and self contained proof of the fact 
the two quotients are homeomorphic,
to us
one of the hardest points in the proof is to show that the mapping
$\chi^{-1}_{\D}$ is continous. This 
has been proved by F. Kirwan in \cite{kirwan}, in the case of
a reductive group acting on a compact K\"ahler manifold with finite
isotropy and in broad generality by F. Loose and P. Heinzner in \cite{hl}.
In the case of a torus acting on a vector space, with homogeneous
moment mapping, the result was proved by A. Neeman
by showing that the zero set of the moment mapping
is a deformation retraction of the vector space \cite{neeman}. 
Based on Neeman's result, G. Schwarz, in his review \cite{schwarz}, 
gave a proof for the general case of reductive groups.
We point out that Neeman's result apply 
directly to our context in the case of the cone $X_{\Sigma_F}$, since the 
proof is based on estimates and computations at the Lie algebra level.
}\end{remark}
\begin{remark}{\rm Observe that the quotient $X_{\D}$, together with
its complex structure, only depends on our choice of generators of the
fan and on the choice of quasilattice, whilst the symplectic
structure depends of course on the polytope, as in the rational case.}
\end{remark}

\end{document}